\documentclass[reqno]{amsart}
\usepackage{amssymb}
\usepackage{amsxtra}
\usepackage[mathscr]{eucal}

\DeclareMathOperator*{\Cat}{\mathbf{Cat}}

\newcommand{\Z}{{\mathbf Z}}
\newcommand{\R}{{\mathbf R}}
\newcommand{\Q}{{\mathbf Q}}
\newcommand{\C}{{\mathbf C}}
\newcommand{\eop}{\hfill$\square$}

\theoremstyle{plain}
\newtheorem{Thm}{Theorem}

\newtheorem{Cor}{Corollary}
\newtheorem{Lem}{Lemma}
\newtheorem{Prop}{Proposition}
\theoremstyle{definition}
\newtheorem{Def}{Definition}
\theoremstyle{remark}
\newtheorem{Rem}{Remark}
\newtheorem*{Notn}{Notation and Terminology}

\numberwithin{equation}{section}

\begin{document}

\title[Multiple $q$-Zeta Values]{Multiple $q$-Zeta Values}

\date{\today}

\author{David~M. Bradley}
\address{Department of Mathematics \& Statistics\\
         University of Maine\\
         5752 Neville Hall
         Orono, Maine 04469-5752\\
         U.S.A.}
\email[]{bradley@math.umaine.edu, dbradley@member.ams.org}

\subjclass{Primary: 11M41; Secondary: 11M06, 05A30, 33D15, 33E20,
30B50}

\keywords{Multiple harmonic series, $q$-analog, multiple zeta
values, $q$-series, Lambert series.}

\begin{abstract} We introduce a $q$-analog 
of the multiple harmonic series commonly referred to as multiple
zeta values.  The multiple $q$-zeta values satisfy a $q$-stuffle
multiplication rule analogous to the stuffle multiplication rule
arising from the series representation of ordinary multiple zeta
values.  Additionally, multiple $q$-zeta values can be viewed as
special values of the multiple $q$-polylogarithm, which admits a
multiple Jackson $q$-integral representation whose limiting case
is the Drinfel'd simplex integral for the ordinary multiple
polylogarithm when $q=1$. The multiple Jackson $q$-integral
representation for multiple $q$-zeta values leads to a second
multiplication rule satisfied by them, referred to as a
$q$-shuffle. Despite this, it appears that many numerical
relations satisfied by ordinary multiple zeta values have no
interesting $q$-extension.  For example, a suitable $q$-analog of
Broadhurst's formula for $\zeta(\{3,1\}^n)$, if one exists, is
likely to be rather complicated.  Nevertheless, we show that a
number of infinite classes of relations, including Hoffman's
partition identities, Ohno's cyclic sum identities, Granville's
sum formula, Euler's convolution formula, Ohno's generalized
duality relation, and the derivation relations of Ihara and Kaneko
extend to multiple $q$-zeta values.
\end{abstract}

\maketitle

\tableofcontents \interdisplaylinepenalty=500

\section{Introduction}\label{sect:Intro}

Throughout, we assume $q$ is real and $0<q<1$.  The $q$-analog of
a non-negative integer $n$ is
\[
   [n]_q := \sum_{k=0}^{n-1} q^k = \frac{1-q^n}{1-q}.
\]
\begin{Def}\label{def:qMZV} Let $m$ be a positive integer and let
$s_1,s_2,\dots,s_m$ be real numbers with $s_1>1$ and $s_j\ge 1$
for $2\le j\le m$.  The multiple $q$-zeta function is the nested
infinite series
\begin{equation}
   \zeta[s_1,\dots,s_m] := \sum_{k_1>\cdots >k_m>0}
   \; \prod_{j=1}^m \frac{q^{(s_j-1)k_j}}{[k_j]_q^{s_j}},
   \label{qMZVdef}
\end{equation}
where the sum is over all positive integers $k_j$ satisfying the
indicated inequalities. If $m=0$, the argument list
in~\eqref{qMZVdef} is empty, and we define $\zeta[\;]:=1$.  If the
arguments in~\eqref{qMZVdef} are positive integers (with $s_1>1$
for convergence), we refer to~\eqref{qMZVdef} as a \emph{multiple
$q$-zeta value.}
\end{Def}

Clearly, $\displaystyle\lim_{q\to 1} \zeta[s_1,\dots,s_m]
=\zeta(s_1,\dots,s_m)$, where
\begin{equation}
   \zeta(s_1,\dots,s_m) := \sum_{k_1>\cdots >k_m>0}
   \; \prod_{j=1}^m k_j^{-s_j},
   \label{MZVdef}
\end{equation}
is the ordinary multiple zeta
function~\cite{BBB,BBBLa,BBBLc,BowBrad1,BowBradSurvey,BowBrad3,BowBradRyoo,Prtn,Hoff92,HoffAlg}.
In this paper, we make a detailed study of the multiple $q$-zeta
function and its values at positive integer arguments.  The
$q$-stuffle rule and some of its implications are worked out in
Section~\ref{sect:qStuffles}.  Among other things, we derive a
$q$-analog of the Newton
recurrence~\cite[eq.~(4.5)]{BowBradSurvey} for $\zeta(\{s\}^n)$, a
$q$-analog of Hoffman's partition identity~\cite[Theorem
2.2]{Hoff92}~\cite{Prtn} and a $q$-analog of the parity reduction
theorem~\cite[Theorem 3.1]{BBBLa}. In Section~\ref{sect:qDuality},
we prove a $q$-analog of Ohno's generalized duality
relation~\cite{Ohno}. Consequences of our generalized $q$-duality
relation include a $q$-analog of ordinary duality for multiple
zeta values, and a $q$-analog of the sum formula~\cite{Granville}.
In Section~\ref{sect:qDerivations}, we prove that the derivation
theorem of Ihara and Kaneko~\cite{Ihara} also extends to multiple
$q$-zeta values.  As we shall see, the $q$-analog of the
Ihara-Kaneko derivation theorem is in fact equivalent to
generalized $q$-duality. A special case ($n=1$) yields a
$q$-analog of Hoffman's derivation relation~\cite[Theorem
5.1]{Hoff92}~\cite[Theorem 2.1]{HoffOhno}. In
Section~\ref{sect:cyclic}, we derive a $q$-analog of Ohno's cyclic
sum formula~\cite{HoffOhno}. In Section~\ref{sect:polylogs}, we
introduce the multiple $q$-polylogarithm, derive a Jackson
$q$-integral analog of the Drinfel'd integral representation for
ordinary multiple polylogarithms, and prove a $q$-analog of a
formula~\cite[Theorem 9.1]{BBBLa} for the colored multiple
polylogarithm. Finally, in Section~\ref{sect:qEuler} we employ
Heine's summation formula for the basic hypergeometric function to
derive a bivariate generating function identity for the multiple
$q$-zeta values $\zeta[m+2,\{1\}^n]$ $(0\le m,n\in\Z)$.  These are
the values of the multiple $q$-zeta function evaluated at the
indecomposable sequences~\cite{Hoff92} consisting of a positive
integer greater than 1 followed by a string of $n$ 1's.
Consequences of our generating function identity include the
special case $\zeta[m+2,\{1\}^n]=\zeta[n+2,\{1\}^m]$ of
$q$-duality, and a $q$-analog of Euler's evaluation expressing
$\zeta(m+2,1)$ as a convolution of ordinary Riemann zeta values.
More generally, we'll see that for all integers $m\ge 2$ and $n\ge
0$, $\zeta[m,\{1\}^n]$ can be expressed in terms of $q$-zeta
values of a single argument. Euler's formula is but a special
case, as is Markett's formula~\cite{Markett} for $\zeta(m,1,1)$.

Whereas the structure of our arguments in many cases derives from
the corresponding arguments in the classical $q=1$ case, the
reader should not be surprised to learn that, as is often the case
with those afflicted with a $q$-virus, much of the difficulty in
establishing an appropriate $q$-theory is determining ``where to
put the $q$.'' In this light, it may be worth remarking that
alternative definitions of the multiple $q$-zeta value are
possible, and lead to other results. For example,
in~\cite{DBqKarl} we study the relationship between certain sums
involving $q$-binomial coefficients with the finite sums
\[
   Z_n[s_1,\dots,s_m] := \sum_{n\ge k_1\ge\cdots\ge k_m\ge 1}\;
   \prod_{j=1}^m \frac{q^{k_j}}{[k_j]_q^{s_j}},
\]
special cases of which have occurred in connection with some
problems in sorting theory.  Another model,
\[
   \zeta_q^{*}(s_1,\dots,s_m) :=
   \sum_{k_1>\cdots>k_m>0}\;\prod_{j=1}^m
   \frac{q^{k_js_j}}{\big(1-q^{k_j}\big)^{s_j}}
\]
is suggested by Zudilin~\cite{Zud}.  See also~\cite{Schles}.  In
Kaneko et al~\cite{Kaneko}, analytic properties of the $q$-analog
\[
   \zeta[s] = \sum_{k=1}^\infty \frac{q^{(s-1)k}}{[k]_q^s}
\]
of the Riemann zeta function are studied.  This immediately
suggested Definition~\ref{def:qMZV} to the present author.
However, as we were subsequently informed, Zhao~\cite{Zhao} had
already been studying~\eqref{qMZVdef} and its polylogarithmic
extension, albeit primarily from the viewpoint of analytic
continuation and the $q$-shuffles of~\cite{BowBradSurvey}. For
arithmetical results on single $q$-zeta values, see
Zudilin~\cite{Zud2,Zud3,Zud4,Zud5,Zud6}.

\begin{Notn}  As customary the boldface symbols $\Z$, $\Q$ and
$\C$ denote the sets of integers, rational numbers, and complex
numbers, respectively.  We'll use $\Z^{+}$ for the set of positive
integers; the subset $\{1,2,\dots,n\}$ consisting of the first $n$
positive integers will be denoted by $\langle n\rangle$. We denote
the cardinality of a set $A$ by $|A|$, and when $A$ is finite, the
group of $|A|!$ permutations of $\langle |A|\rangle$ by
${\mathfrak{S}(A)}$.  If $A=\langle n\rangle$, we write
${\mathfrak{S}}_n$ instead of $\mathfrak{S}(\langle n\rangle)$.
Boolean expressions such as $(k\in A)$ take the value $1$ if $k\in
A$ and $0$ if $k\notin A$. To avoid the potential for ambiguity in
expressing complicated argument sequences without recourse to
ellipses, we make occasional use of the abbreviations
$\Cat_{j=1}^m \{s_j\}$ for the concatenated argument sequence
$s_1,\dots,s_m$ and $\{s\}^m=\Cat_{j=1}^m \{s\}$ for $m\ge 0$
consecutive copies of $s$, which may itself be a sequence of
arguments.  Throughout, $I$ will denote the set $\{0,1\}$ and
$I^m$ the Cartesian product $I\times\cdots \times I$ of $m$ copies
of $I$ when $m$ is a positive integer. This will cause no
confusion with the notation for concatenation, since we will never
have occasion to discuss the periodic sequence $0,1,\dots,0,1$. As
in~\cite{BBBLa}, we define the \emph{depth} of the multiple
$q$-zeta function~\eqref{qMZVdef} to be the number $m$ of
arguments.
\end{Notn}

\section{$q$-Stuffles}\label{sect:qStuffles}
The stuffle multiplication rule~\cite{BBBLa,BowBradSurvey,Prtn}
for the multiple zeta function (also referred to as the harmonic
product or $*$-product in~\cite{HoffAlg,HoffOhno}) arises when one
expands the product of two nested series of the
form~\eqref{MZVdef}, and is invariably given a recursive
description.  We begin with an explicit formula for the
$q$-stuffle multiplication rule satisfied by the multiple $q$-zeta
function; an explicit formula for the stuffle rule can then be
derived by taking the limit as $q\to 1$.

Let $m$ and $n$ be positive integers.  Define a \emph{stuffle} on
$(m,n)$ as a pair $(\phi,\psi)$ of order-preserving injective
mappings $\phi:\langle m\rangle \to \langle m+n \rangle $, $\psi:
\langle n\rangle \to \langle m+n\rangle$ such that the union of
their images is equal to $\langle r\rangle$ for some positive
integer $r$ with $\max(m,n)\le r\le m+n$.  In what follows we'll
abuse notation by writing (for example) $\phi^{-1}(k)$ for the
pre-image $\phi^{-1}(\{k\})$ of the singleton $\{k\}$.  Since
$\phi$ is injective, $\phi^{-1}(k)$ is either empty $\{\}$ or a
singleton $\{j\}$ for some positive integer $j$, and we make the
conventions
\[
   s_{\{j\}}= s_j,\quad t_{\{j\}}=t_j,\quad s_{\{\}}=t_{\{\}}=0.
\]
The stuffle multiplication rule for the multiple zeta function can
now be written in the form
\begin{equation}
   \zeta(s_1,\dots,s_m)\zeta(t_1,\dots,t_n)
   = \sum_{(\phi,\psi)}\zeta\big(\Cat_{k=1}^r
   \big\{s_{\phi^{-1}(k)}+t_{\psi^{-1}(k)}\big\}\big),
\label{stuffle}
\end{equation}
where the sum is over all stuffles $(\phi,\psi)$ on $(m,n)$, and
$r=r(\phi,\psi)$ is the cardinality (equivalently, the largest
member) of the union $\phi(\langle m\rangle)\cup\psi(\langle
n\rangle)$ of the images of $\phi$ and $\psi$. More generally,
expanding the product
\[
   \zeta[s_1,\dots,s_m]\zeta[t_1,\dots,t_n]
   = \sum_{k_1>\cdots >k_m>0}\;\prod_{j=1}^m
   \frac{q^{(s_j-1)k_j}}{[k_j]_q^{s_j}}
   \sum_{l_1>\cdots >l_n>0}\;\prod_{j=1}^n
   \frac{q^{(t_j-1)l_j}}{[l_j]_q^{t_j}}
\]
yields sums of products of terms of the form
\[
   \frac{q^{(s-1)k+(t-1)l}}{[k]_q^s[l]_q^t},
\]
which, if $k=l$, reduces to
\[
   \frac{q^{(s+t-2)k}}{[k]_q^{s+t}}
   =
   (1-q)\frac{q^{(s+t-2)k}}{[k]_q^{s+t-1}}+\frac{q^{(s+t-1)k}}{[k]_q^{s+t}}.
\]
It follows that
\begin{equation}
   \zeta[s_1,\dots,s_m]\zeta[t_1,\dots,t_n]
   = \sum_{(\phi,\psi)}\sum_{A}
   (1-q)^{|A|}\, \zeta\big[\Cat_{k=1}^r \big\{
   s_{\phi^{-1}(k)}+t_{\psi^{-1}(k)}-(k\in A)\big\}\big],
\label{qstuffle}
\end{equation}
where the outer sum is over all stuffles $(\phi,\psi)$ on $(m,n)$,
the inner sum is over all subsets $A$ of the intersection of the
images of $\phi$ and $\psi$, $r=|\phi(\langle
m\rangle)\cup\psi(\langle n\rangle)|$ as in~\eqref{stuffle}, and
the Boolean expression $(k\in A)$ takes the value $1$ if $k\in A$
and $0$ if $k\notin A$. We refer to~\eqref{qstuffle} as the
$q$-stuffle multiplication rule. Note that~\eqref{stuffle} is the
limiting case $q\to 1$ of~\eqref{qstuffle}.  For an alternative
$q$-deformation of the stuffle algebra, see~\cite{HoffComb}.

\subsection{Period-1 Sums Completely Reduce}
As an application of the $q$-stuffle multiplication
rule~\eqref{qstuffle}, we show that for any $s>1$ and positive
integer $n$, the multiple $q$-zeta function $\zeta[\{s\}^n]$ can
be expressed polynomially in terms of $q$-zeta functions of depth
1.  See~\cite{BowBrad1} for a discussion of the period-2 case for
ordinary multiple zeta values and related alternating Euler sums.
\begin{Thm}\label{thm:zsn}
If $n$ is a positive integer and $s>1$, then
\[
   n\zeta[\{s\}^n] = \sum_{k=1}^n (-1)^{k+1}\zeta[\{s\}^{n-k}]
   \sum_{j=0}^{k-1}\binom{k-1}{j}(1-q)^j\zeta[ks-j].
\]
\end{Thm}
\noindent{\bf Proof.} Let $R$ denote the right-hand side of the
equation in Theorem~\ref{thm:zsn}.   The $q$-stuffle
multiplication rule~\eqref{qstuffle} implies that
\begin{multline}
   R = \sum_{k=1}^n (-1)^{k+1}\sum_{j=0}^{k-1}\binom{k-1}{j}(1-q)^j
   \bigg\{\sum_{m=0}^{n-k}\zeta[\{s\}^m,ks-j,\{s\}^{n-k-m}]\\
   +\sum_{m=0}^{n-1-k}\zeta[\{s\}^m,(k+1)s-j,\{s\}^{n-1-k-m}]\\
   +(1-q)\sum_{m=0}^{n-1-k}\zeta[\{s\}^m,(k+1)s-j-1,\{s\}^{n-1-k-m}]\bigg\}.
\label{R3a}
\end{multline}
Now expand~\eqref{R3a} into three triple sums.  We re-index the
first and third of these, replacing $k$ by $k+1$ in the first, and
$j$ by $j-1$ in the third.  Then
\begin{equation}\label{R3b}
\begin{split}
   R &= \sum_{k=0}^{n-1}(-1)^k \sum_{j=0}^k \binom{k}{j}(1-q)^j
   \sum_{m=0}^{n-1-k}\zeta[\{s\}^m,(k+1)s-j,\{s\}^{n-1-k-m}]\\
   &+\sum_{k=1}^{n-1} (-1)^{k+1}\sum_{j=0}^{k-1}\binom{k-1}{j}(1-q)^j
   \sum_{m=0}^{n-1-k}\zeta[\{s\}^m,(k+1)s-j,\{s\}^{n-1-k-m}]\\
   &+\sum_{k=1}^{n-1}
   (-1)^{k+1}\sum_{j=1}^k\binom{k-1}{j-1}(1-q)^j\sum_{m=0}^{n-1-k}
   \zeta[\{s\}^m,(k+1)s-j,\{s\}^{n-1-k-m}].
\end{split}
\end{equation}
In the second and third triple sums~\eqref{R3b}, we have omitted
the terms corresponding to $k=n$, because these vanish.  In the
second triple sum~\eqref{R3b}, the range on $j$ can be extended to
include the term $j=k$ because the binomial coefficient vanishes
in that case. Similarly, the range on $j$ in the third
sum~\eqref{R3b} can be extended to include the term $j=0$.  If we
now combine the extended second and third triple sums~\eqref{R3b}
using the Pascal formula
\[
   \binom{k-1}{j}+\binom{k-1}{j-1}=\binom{k}{j},
\]
we see that
\begin{equation}\label{R3c}
\begin{split}
   R &= \sum_{k=0}^{n-1}(-1)^k \sum_{j=0}^k \binom{k}{j}(1-q)^j
   \sum_{m=0}^{n-1-k}\zeta[\{s\}^m,(k+1)s-j,\{s\}^{n-1-k-m}]\\
   &+\sum_{k=1}^{n-1} (-1)^{k+1}\sum_{j=0}^{k}\binom{k}{j}(1-q)^j
   \sum_{m=0}^{n-1-k}\zeta[\{s\}^m,(k+1)s-j,\{s\}^{n-1-k-m}].
\end{split}
\end{equation}
The two triple sums~\eqref{R3c} cancel except for the $k=0$ term
in the first.  Thus, we find that
\[
   R=\sum_{m=0}^{n-1}\zeta[\{s\}^m,s,\{s\}^{n-1-m}]
   = n\zeta[\{s\}^n],
\]
as required. \eop

For reference, we note that letting $q\to 1$ in
Theorem~\ref{thm:zsn} yields the Newton recurrence~\cite[eq.\
(4.5)]{BowBradSurvey} for multiple zeta values of period 1.
\begin{Cor}
If $n$ is a positive integer and $s>1$, then
\[
   n\zeta(\{s\}^n) = \sum_{k=1}^n (-1)^{k+1}\zeta(\{s\}^{n-k})
   \zeta(ks).
\]
\end{Cor}

\subsection{Partition Identities}

Additional $q$-stuffle relations can be most easily stated using
the concept of a set partition.  As in~\cite{Prtn}, it is helpful
to distinguish between set partitions that are ordered and those
that are unordered.

\begin{Def}[Unordered Set Partition]
Let $S$ be a finite non-empty set.  An \emph{unordered} set
partition of $S$ is a finite non-empty set $\mathscr{P}$ whose
elements are disjoint non-empty subsets of $S$ with union $S$.
That is, there exists a positive integer $m=|\mathscr{P}|$ and
non-empty subsets $P_1,\dots,P_m$ of $S$ such that
$\mathscr{P}=\{P_1,\dots,P_m\}$, $S=\cup_{k=1}^m P_k$, and
$P_j\cap P_k$ is empty if $j\ne k$.
\end{Def}

\begin{Def}[Ordered Set Partition]
Let $S$ be a finite non-empty set.  An \emph{ordered} set
partition of $S$ is a finite ordered tuple $\vec P$ of disjoint
non-empty subsets of $S$ such that the union of the components of
$\vec P$ is equal to $S$. That is, there exists a positive integer
$m$ and non-empty subsets $P_1,\dots,P_m$ of $S$ such that $\vec
P$ can be identified with the ordered $m$-tuple $(P_1,\dots,P_m)$,
$\cup_{k=1}^m P_k=S$, and $P_j\cap P_k$ is empty if $j\ne k$.
\end{Def}

We next introduce the shift operators $E_k$ and $\delta_k$ defined
as follows.
\begin{Def} Let $m$ and $k$ be positive integers with $1\le k\le
m$, and let $s_1,\dots,s_m$ be real numbers with $s_1>1$, $s_k\ge
2$ and $s_j\ge 1$ for $2\le j\ne k\le m$. The shift operator $E_k$
is defined by means of
\[
   E_k \zeta[s_1,\dots,s_m] = \zeta\big[\Cat_{j=1}^{k-1}
   s_j,s_k-1,\Cat_{j=k+1}^m s_j\big].
\]
Let $\delta_k:= \delta_k(q)=1+(1-q)E_k$ and abbreviate
$\delta:=\delta_1$.
\end{Def}
The $q$-stuffle multiplication rule~\eqref{qstuffle} can now be
re-written in the form
\begin{equation}
   \zeta[s_1,\dots,s_m]\zeta[t_1,\dots,t_n]
   = \sum_{(\phi,\psi)}
   \bigg(\prod_{k=1}^r \delta_k^{\alpha_k}\bigg)
   \zeta\big[\Cat_{k=1}^r
   \big\{s_{\phi^{-1}(k)}+t_{\psi^{-1}(k)}\big\}\big],
\label{deltaqstuffle}
\end{equation}
where $r=|\phi(\langle m\rangle)\cup\psi(\langle n\rangle)|$ and
$\alpha_k$ is equal to 1 or 0 according as to whether $k$
respectively is or is not a member of the intersection
$\phi(\langle m\rangle)\cap\psi(\langle n\rangle)$ of the images
of $\phi$ and $\psi$.  Given~\eqref{deltaqstuffle}, the following
result is self-evident, but it can also be readily proved by
mathematical induction.
\begin{Thm}\label{thm:nproduct}
Let $n$ be a positive integer, and let $s_k>1$ for
$1\le k\le n$.  Then
\begin{align*}
   \prod_{k=1}^n \zeta[s_k] &=
   \sum_{\vec P\vDash \langle n\rangle}
   \bigg(\prod_{j=1}^{|\vec P|}\delta_j^{|P_j|-1}\bigg)
   \zeta\bigg[\Cat_{j=1}^{|\vec P|}
   \sum_{i\in P_j}s_i\bigg]\\
   &= \sum_{m=1}^n
   \sum_{\substack{\vec P\vDash \langle n\rangle\\ |\vec P|=m}}
   \sum_{\nu_1=0}^{|P_1|-1}\cdots\sum_{\nu_m=0}^{|P_m|-1}
   \zeta\bigg[\Cat_{j=1}^m \sum_{i\in P_j}s_i-\nu_j\bigg] \prod_{j=1}^m
   \binom{|P_j|-1}{\nu_j}(1-q)^{\nu_j},
\end{align*}
where the sum is over all ordered set partitions $\vec P$ of
$\langle n\rangle$ having components $(P_1,\dots,P_m)$, with $1\le
m=|\vec P|\le n$.
\end{Thm}
If in Theorem~\ref{thm:nproduct} we abbreviate $\sum_{i\in P_j}
s_i$ by $p_j$ and sum instead over unordered set partitions, we
see that
\begin{equation}
   \prod_{k=1}^n \zeta[s_k] =
    \sum_{\mathscr{P}\vdash \langle n\rangle}
   \bigg(\prod_{j=1}^{|\mathscr{P}|} \delta_j^{|P_j|-1}\bigg)
   \sum_{\sigma\in{\mathfrak{S}(\mathscr{P})}}
   \zeta\bigg[\Cat_{j=1}^{|\mathscr{P}|} p_{\sigma(j)}\bigg],
\label{nproductperm}
\end{equation}
where the $P_j\subseteq\langle n\rangle$ are the distinct disjoint
members of $\mathscr{P}$.  Inverting~\eqref{nproductperm} and
expanding the delta operators yields the following partition
identity.
\begin{Thm}\label{thm:qHoffman} Let $n$ be a positive integer, and
let $s_j>1$ for $1\le j\le n$.  Then
\[
   \sum_{\sigma\in{\mathfrak{S}_n}} \zeta\big[\Cat_{j=1}^n
   s_{\sigma(j)}\big]= \sum_{\mathscr{P}\vdash \langle n\rangle}
   (-1)^{n-|\mathscr{P}|} \prod_{k=1}^{|\mathscr{P}|} (|P_k|-1)!
   \sum_{\nu_k=0}^{|P_k|-1}\binom{|P_k|-1}{\nu_k}(1-q)^{\nu_k}
   \zeta[p_k-\nu_k],
\]
where the sum on the right is over all unordered set partitions
$\mathscr{P}=\{P_1,\dots,P_m\}$ of $\langle n\rangle$, $1\le
m=|\mathscr{P}|\le n$, and $p_k= \sum_{j\in P_k}s_j$.
\end{Thm}

Letting $q\to 1$ in Theorem~\ref{thm:qHoffman}, we obtain the
following result of Hoffman~\cite[Theorem 2.2]{Hoff92}, which he
proved using a counting argument.
\begin{Cor}[Hoffman's Partition Identity]
\label{cor:Hoffman} Let $n$ be a positive integer, and let $s_j>1$
for $1\le j\le n$.  Then
\[
   \sum_{\sigma\in{\mathfrak{S}_n}} \zeta\big(\Cat_{j=1}^n
   s_{\sigma(j)})= \sum_{\mathscr{P}\vdash \langle n\rangle}
   (-1)^{n-|\mathscr{P}|} \prod_{P\in\mathscr{P}} (|P|-1)!\,
   \zeta\big(\sum_{j\in P}s_j\big),
\]
where the sum on the right is over all unordered set partitions
$\mathscr{P}$ of $\langle n\rangle$.
\end{Cor}

\noindent{\bf Proof of Theorem~\ref{thm:qHoffman}.} It is enough
to show that
\begin{equation}
   \sum_{\sigma\in{\mathfrak{S}_n}} \zeta\big[\Cat_{j=1}^n
   s_{\sigma(j)}]= \sum_{\mathscr{P}\vdash \langle n\rangle}
   (-1)^{n-|\mathscr{P}|} \prod_{P\in\mathscr{P}} (|P|-1)!\,
   \delta^{|P|-1}\zeta\bigg[\sum_{j\in P}s_j\bigg].
\label{qHoffman}
\end{equation}
When $n=1$ this is trivial.  Suppose the result~\eqref{qHoffman}
holds for $n-1$.  Then
\begin{equation}
   \sum_{\sigma\in{\mathfrak{S}_{n-1}}} \zeta\big[\Cat_{j=1}^{n-1}
   s_{\sigma(j)}]= \sum_{\mathscr{P}\vdash \langle n-1\rangle}
   (-1)^{n-1-|\mathscr{P}|} \prod_{P\in\mathscr{P}} (|P|-1)!\,
   \delta^{|P|-1}\zeta\bigg[\sum_{j\in P}s_j\bigg].
\label{qHoffindhyp}
\end{equation}
After multiplying equation~\eqref{qHoffindhyp} through by
$\zeta[s_n]$, applying the $q$-stuffle multiplication
rule~\eqref{deltaqstuffle} to the left hand side, and moving the
stuffed terms to the right, we obtain
\begin{multline}
    \sum_{\sigma\in{\mathfrak{S}_n}} \zeta\big[\Cat_{j=1}^n
   s_{\sigma(j)}]=\sum_{\mathscr{P}\vdash \langle n-1\rangle}
   (-1)^{n-1-|\mathscr{P}|} \zeta[s_n]\prod_{P\in\mathscr{P}} (|P|-1)!\,
   \delta^{|P|-1}\zeta\bigg[\sum_{j\in P}s_j\bigg]\\
   -\sum_{\sigma\in{\mathfrak{S}_{n-1}}}\sum_{k=1}^{n-1} \delta_k
   \zeta\big[\Cat_{j=1}^{k-1}s_{\sigma(j)},s_{\sigma(j)}+s_n,
   \Cat_{j=k+1}^{n-1}s_{\sigma(j)}\big].
\label{bigmess}
\end{multline}
Let $u_j^{(k)}=s_j$ if $j\ne k$ and $u_k^{(k)}=s_k+s_n$.  With the
aid of the inductive hypothesis~\eqref{qHoffindhyp}, the double
sum on the right hand side of~\eqref{bigmess} can now be expressed
in the form
\[
   \sum_{k=1}^{n-1}\sum_{\sigma\in{\mathfrak{S}_n}}
   \delta_{\sigma^{-1}(k)}\zeta[\Cat_{j=1}^{n-1}u_{\sigma(j)}^{(k)}]
   = \sum_{k=1}^{n-1}\sum_{\mathscr{P}\vdash\langle n-1\rangle}
   (-1)^{n-1-|\mathscr{P}|}\prod_{P\in\mathscr{P}}(|P|-1)!\,
   \delta^{|P|-1+(k\in P)}\zeta\big[\sum_{j\in P}u_j^{(k)}\big].
\]
From~\eqref{bigmess}, it now follows that
\begin{multline}
    \sum_{\sigma\in{\mathfrak{S}_n}} \zeta\big[\Cat_{j=1}^n
   s_{\sigma(j)}]=\sum_{\mathscr{P}\vdash \langle n-1\rangle}
   (-1)^{n-1-|\mathscr{P}|} \zeta[s_n]\prod_{P\in\mathscr{P}} (|P|-1)!\,
   \delta^{|P|-1}\zeta\bigg[\sum_{j\in P}s_j\bigg]\\
   +\sum_{k=1}^{n-1}\sum_{\mathscr{P}\vdash\langle n-1\rangle}
   (-1)^{n-|\mathscr{P}|}\prod_{P\in\mathscr{P}}(|P|-1)!\,
   \delta^{|P|-1+(k\in P)}\zeta\bigg[\sum_{j\in P}u_j^{(k)}\bigg].
\label{littlemess}
\end{multline}
Note that in the second sum on the right hand side
of~\eqref{littlemess}, as $k$ runs from 1 to $n-1$, there is a
contribution of $|P_0|$ copies of the inner sum if $P_0\in
\mathscr{P}$ is such that $k\in P_0$.  Therefore, if to each
partition $\mathscr{P}$ of $\langle n-1\rangle$ in the first sum
on the right hand side of~\eqref{littlemess}, we let
$\mathscr{R}=\mathscr{P}\cup\{\{n\}\}$, then
\begin{multline}
    \sum_{\sigma\in{\mathfrak{S}_n}} \zeta\big[\Cat_{j=1}^n
   s_{\sigma(j)}]=\sum_{\substack{\mathscr{R}\vdash \langle
   n\rangle\\ \{n\}\in \mathscr{R}}}
   (-1)^{n-|\mathscr{R}|} \prod_{R\in\mathscr{R}} (|R|-1)!\,
   \delta^{|R|-1}\zeta\bigg[\sum_{j\in R}s_j\bigg]\\
   \sum_{\substack{\mathscr{P}\vdash\langle n-1\rangle\\ P_0\in\mathscr{P}}}
   (-1)^{n-|\mathscr{P}|}
   |P_0|!\,\delta^{|P_0|}\zeta\big[s_n+\sum_{j\in
   P_0}s_j\big]
   \prod_{\substack{P\in\mathscr{P}\\P\ne P_0}}(|P|-1)!\,
   \delta^{|P|-1}\zeta\bigg[\sum_{j\in P}s_j\bigg].
\label{mess}
\end{multline}
Clearly, the second sum on the right hand side of~\eqref{mess} can
be re-written more succinctly if we simply toss $n$ into $P_0$ and
thus view each $\mathscr{P}$ as an unordered set partition of
$\langle n\rangle$ in which no part in the partition is equal to
the singleton $\{n\}$.  Thus,
\begin{multline*}
    \sum_{\sigma\in{\mathfrak{S}_n}} \zeta\big[\Cat_{j=1}^n
   s_{\sigma(j)}]
   =\sum_{\substack{\mathscr{R}\vdash \langle
   n\rangle\\ \{n\}\in \mathscr{R}}}
   (-1)^{n-|\mathscr{R}|} \prod_{R\in\mathscr{R}} (|R|-1)!\,
   \delta^{|R|-1}\zeta\bigg[\sum_{j\in R}s_j\bigg]\\
   \sum_{\substack{\mathscr{P}\vdash\langle n\rangle\\
   \{n\}\notin\mathscr{P}}}(-1)^{n-|\mathscr{P}|}
   \prod_{P\in\mathscr{P}}(|P|-1)!\,
   \delta^{|P|-1}\zeta\bigg[\sum_{j\in P}s_j\bigg].
\end{multline*}
The result~\eqref{qHoffman} now follows, since any partition of
$\langle n \rangle$ is either of the form $\mathscr{R}$ or
$\mathscr{P}$ above.   \eop

\begin{Rem}
The proof shows that Theorem~\ref{thm:qHoffman} (and hence also
its limiting case, Corollary~\ref{cor:Hoffman}) relies on only the
$q$-stuffle multiplication property.  Loosely speaking, we refer
to results such as Theorems~\ref{thm:nproduct}
and~\ref{thm:qHoffman} and Corollary~\ref{cor:Hoffman} as
\emph{partition identities} because they are easily stated using
the language of set partitions.  The notion is defined precisely
in~\cite{Prtn}, where among other things it is shown that
\emph{all} partition identities are a consequence of the stuffle
multiplication rule, and hence a decision procedure exists for
verifying them.
\end{Rem}

We conclude this section with one further result, namely a
$q$-analog of~\cite[Theorem 3.1]{BBBLa}. Results which go beyond
stuffles will be discussed in the subsequent sections.

\begin{Thm}[Parity Reduction]
Let $m$ and let $s_1,\dots,s_m$ be real numbers with $s_1>1$,
$s_m>1$, and $s_j\ge 1$ for $1<j<m$. Then
\[
  \zeta\big[\Cat_{k=1}^m s_k\big]
  +(-1)^m \zeta\big[\Cat_{k=1}^m s_{m-k+1}\big]
\]
can be expressed as a $\Z[q]$-linear combination of multiple
$q$-zeta values of depth less than $m$.  That is, the coefficients
in the linear combination are polynomials in $q$ with integer
coefficients.
\end{Thm}
\noindent{\bf Proof.} Let $N$ denote the Cartesian product of $m$
copies of the positive integers.  Define an additive
weight-function on subsets of $N$ by
\[
   w(A) := \sum_{\vec n\in A}\; \prod_{k=1}^m \frac{q^{(s_k-1)n_k}}
   {[n_k]_q^{s_k}},
\]
where the sum is over all $\vec n=(n_1,\dots,n_m)\in A$.   For
each $k\in\langle m-1\rangle$, define the subset $P_k$ of $N$ by
$P_k=\{\vec n\in N: n_k\le n_{k+1}\}$. The Inclusion-Exclusion
Principle states that
\begin{equation}
\label{inc}
   w\bigg(\bigcap_{k=1}^{m-1} N\setminus P_k\bigg)
   = \sum_{T\subseteq\langle m-1\rangle}(-1)^{|T|}\;
   w\bigg(\bigcap_{k\in T} P_k\bigg).
\end{equation}
The term on the right hand side of~\eqref{inc} arising from the
empty subset $T=\{\}$ is $\prod_{k=1}^m \zeta[s_k]$ by the usual
convention for intersection over an empty set.  The left hand side
of~\eqref{inc} is simply $\zeta[s_1,\dots,s_m]$.  In light of the
identity
\[
   \frac{q^{(s-2)n}}{[n]_q^s}=\frac{q^{(s-1)n}}{[n]_q^s}+
   (1-q)\frac{q^{(s-2)n}}{[n]_q^{s-1}},
\]
it follow that all terms on the right hand side of~\eqref{inc} are
multiple $q$-zeta values of depth strictly less than $m$, except
when $T=\langle m-1\rangle$, which contributes
\begin{multline*}
   (-1)^{m-1}\sum_{1\le n_1\le n_2\le\cdots\le n_m}\;
   \prod_{k=1}^m \frac{q^{(s_k-1)n_k}}{[n_k]_q^{s_k}}\\
   =(-1)^{m-1}\zeta\big[\Cat_{k=1}^m
   s_{m-k+1}\big]+ \mbox{lower depth multiple $q$-zeta values}.
\end{multline*}
\eop

\section{Generalized $q$-Duality}\label{sect:qDuality}
In this section, we prove a $q$-analog of Ohno's generalized
duality relation~\cite{Ohno}.  As a consequence, we derive
$q$-analogs of the duality relation and the sum
formula~\cite{Granville}.  An additional consequence is a
$q$-analog of Ihara and Kaneko's derivation theorem~\cite{Ihara},
which we prove in Section~\ref{sect:qDerivations}.

\begin{Def} Let $n$ and $s_1,\dots,s_n$ be positive integers with
$s_1>1$.  Let $m$ be a non-negative integer.  Define
\[
   Z[s_1,\dots,s_n;m] := \sum_{\substack{c_1,\dots,c_n\ge
   0\\c_1+\cdots+c_n=m}} \zeta[s_1+c_1,\dots,s_n+c_n],
\]
where the sum is over all non-negative integers $c_j$ with
$\sum_{j=1}^n c_j=m$.  As in~\cite{BBB}, for non-negative integers
$a_j$ and $b_j$, define the dual argument lists
\[
   \mathbf p = \big(\Cat_{j=1}^n \{a_j+2,\{1\}^{b_j}\}\big),
   \qquad
   \mathbf p' =\big(\Cat_{j=1}^n
   \{b_{n-j+1}+2,\{1\}^{a_{n-j+1}}\}\big).
\]
\end{Def}

\begin{Thm}[Generalized $q$-duality]\label{thm:genqduality}
For any pair of dual argument lists $\mathbf p$, $\mathbf p'$ and
any non-negative integer $m$, we have the equality $Z[\mathbf
p;m]=Z[\mathbf p';m]$.
\end{Thm}

The $m=0$ case of Theorem~\ref{thm:genqduality} is worth stating
separately.  It is a direct $q$-analog of the duality relation for
multiple zeta values.
\begin{Cor}[$q$-duality]\label{cor:qduality}
For any pair of dual argument lists $\mathbf p$ and $\mathbf p'$,
we have the equality $\zeta[\mathbf p] = \zeta[\mathbf p']$.  In
other words, for all non-negative integers $a_j,b_j$, $1\le j\le
n$, we have the equality
\[
   \zeta\big[\Cat_{j=1}^n \{a_j+2,\{1\}^{b_j}\}\big]
   = \zeta\big[\Cat_{j=1}^n
   \{b_{n-j+1}+2,\{1\}^{a_{n-j+1}}\}\big].
\]
\end{Cor}

As noted by Ohno~\cite{Ohno}, the sum formula~\cite{Granville} is
an easy consequence of his generalized duality relation. Likewise,
the following $q$-analog of the sum formula is a consequence of
our generalized $q$-duality relation
(Theorem~\ref{thm:genqduality}).
\begin{Cor}[$q$-sum formula]\label{cor:qsum} For any integers $0<k\le n$,
we have
\[
   \sum_{s_1+s_2+\cdots+s_n=k}\zeta[s_1+1,s_2,\dots,s_n]
   = \zeta[k+1],
\]
where the sum is over all positive integers $s_1,s_2,\dots,s_n$
with sum equal to $k$.
\end{Cor}
\noindent{\bf Proof.}  If we take the dual argument lists in the
form $\mathbf p = (n+1)$ and $\mathbf p' = (2,\{1\}^{n-1})$ and
put $m=k-n$, then Theorem~\ref{thm:genqduality} states that
\[
   \zeta[k+1]=
   \sum_{\substack{c_1,\dots,c_n\ge 0\\c_1+\cdots+c_n=k-n}}
   \zeta\big[2+c_2,\Cat_{j=2}^n\{1+c_j\}\big]
   =\sum_{\substack{s_1,\dots,s_n\ge 1\\s_1+\cdots+s_n=k}}
   \zeta\big[s_1+1,\Cat_{j=2}^n s_j\big].
\]
\eop

\subsection{Proof of Generalized $q$-Duality}

To prove Theorem~\ref{thm:genqduality}, we need to employ some
algebraic machinery first introduced by Hoffman~\cite{HoffAlg}.
The argument itself extends ideas of Okuda and
Ueno~\cite{OkudaUeno} to the $q$-case.   Let
$\mathfrak{h}=\Q\langle x,y\rangle$ denote the non-commutative
polynomial algebra over the rational numbers in two indeterminates
$x$ and $y$, and let $\mathfrak{h}^0$ denote the subalgebra
$\Q1\oplus x\mathfrak{h}y$. The $\Q$-linear map
$\widehat{\zeta}:\mathfrak{h}^0\to \R$ is defined by
$\widehat{\zeta}[1] := \zeta[\,]=1$ and
\[
   \widehat{\zeta}\bigg[\prod_{j=1}^s x^{a_j}y^{b_j}\bigg]
   =\zeta\big[\Cat_{j=1}^s \{a_j+1,\{1\}^{b_j-1}\}\big],
   \qquad a_j,b_j\in\Z^{+}.
\]
For each positive integer $n$, let $D_n$ be the derivation on
$\mathfrak{h}$ that maps $x\mapsto 0$ and $y\mapsto x^n y$, and
let $\theta$ be a formal parameter.  Then $\sum_{n=1}^\infty D_n
\theta^n/n$ is a derivation on $\mathfrak{h}[[\theta]]$ and
$\sigma_{\theta} = \exp\big(\sum_{n=1}^\infty D_n\theta^n/n\big)$
is an automorphism of $\mathfrak{h}[[\theta]]$.   Let $\tau$ be
the anti-automorphism of $\mathfrak{h}$ that switches $x$ and $y$.
For any word $w\in\mathfrak{h}^0$, define
$f[w;\theta]:=\widehat{\zeta}[\sigma_{\theta}(w)]$ and
$g[w;\theta] := \widehat{\zeta}[\sigma_{\theta}(\tau(w))]$.  By
definition of $D_n$, $\sum_{n=1}^\infty D_n\theta^n/n$ sends
$x\mapsto 0$ and $y\mapsto \{\log(1-x\theta)^{-1}\}y$.  Thus,
$\sigma_{\theta}$ sends $x\mapsto x$ and $y\mapsto
(1-x\theta)^{-1}y$.  Therefore,
\begin{equation}\label{fgf}
   f\bigg[\prod_{j=1}^s x^{a_j}y^{b_j};\theta\bigg]
   = \widehat{\zeta}\bigg[\prod_{j=1}^s
   x^{a_j}\big\{(1-x\theta)^{-1}y\big\}^{b_j}\bigg]
   =\sum_{m=0}^\infty \theta^m\sum_{\substack{c_1,\dots,c_n\ge
   0\\ c_1+\cdots+c_n=m}}\zeta\big[\Cat_{i=1}^n \{k_i+c_i\}\big],
\end{equation}
where $(k_1,\dots,k_n)=(\Cat_{j=1}^s\{a_j+1,\{1\}^{b_j-1}\})$ and
$n=\sum_{j=1}^s b_j$.  Theorem~\ref{thm:genqduality} can now be
restated in the equivalent form given below.
\begin{Thm}[Generalized $q$-duality, reformulated]\label{thm:f=g} For all
$w\in\mathfrak{h}^0$, $f[w;\theta]=g[w;\theta]$.  In other words,
$\widehat{\zeta}\circ\sigma_{\theta}$ is invariant under ordinary
duality $\tau$.
\end{Thm}

The following difference equation is the key result we need to
prove Theorem~\ref{thm:f=g}.
\begin{Thm}\label{thm:qdiff} Let $a_i,b_i$
be positive integers with $\sum_{i=1}^s (a_i+b_i)>2$.  Make the
abbreviation $\theta\,':=q\theta-1$, and recall the notation
$I^m=\{0,1\}\times \cdots\times\{0,1\}$ for the $m$-fold Cartesian
product from Section~\ref{sect:Intro}. The generating functions
$f$ and $g$ satisfy the difference equation
\begin{multline*}
   \sum_{\substack{\epsilon,\delta\in I^s\\ \delta_1<a_1, \epsilon_s<b_s}}
   (-\theta)^{\overline{\delta}\cdot\overline{\epsilon}}(1-q)^{\delta\cdot\epsilon}
   f\bigg[\prod_{i=1}^s
   x^{a_i-\delta_i}y^{b_i-\epsilon_i};\theta\bigg]\\
   =
   \sum_{\substack{\delta,\epsilon\in I^{s+1}\\ \delta_{s+1}=\epsilon_1=0\\ \delta_1<a_1, \epsilon_{s+1}<b_s}}
   (-\theta\,')^{\overline{\delta}\cdot\overline{\epsilon}-1}(1-q)^{\delta\cdot\epsilon}
   f\bigg[\prod_{i=1}^s
   x^{a_i-\delta_i}y^{b_i-\epsilon_{i+1}};\theta\,'\bigg].
\end{multline*}
Here, we use $\overline{\delta}$ to denote the ordered tuple whose
$i^{\mathrm{th}}$ component is $1-\delta_i$, and of course
$\delta\cdot\epsilon$ denotes the dot product $\sum_i
\delta_i\epsilon_i$.  Similarly, $\overline{\epsilon}$ denotes the
ordered tuple whose $i^{\mathrm{th}}$ component is $1-\epsilon$,
and
$\overline{\delta}\cdot\overline{\epsilon}=\sum_{i}(1-\delta_i)(1-\epsilon_i)$.
\end{Thm}

We also require the following lemma, which shows that the
generating function $f[w;\theta]$ can be analytically continued to
a meromorphic function of $\theta$ with at worst simple poles at
$\theta = q^{-\nu}[\nu]_q$ for positive integers $\nu$.
\begin{Lem}\label{lem:fmero} Let $w=\prod_{i=1}^s x^{a_i}y^{b_i}$, where $a_i$ and
$b_i$ are positive integers.  Let $B_0 := 0$ and set $B_i :=
\sum_{j=1}^i b_j$ for $1\le i\le s$.   Then
\[
   f[w;\theta] = \sum_{\nu=1}^\infty \frac{C_{\nu}[w]}{[\nu]_q-\theta
   q^{\nu}},
\]
where
\[
   C_{\nu}[w] := \sum_{k=1}^{B_s}
   \sum_{\substack{m_1>\cdots>m_{k-1}>\nu\\ \nu>m_{k+1}>\cdots>m_{B_s}>0}}
   E_k[w;m_1,\dots,m_{k-1},\nu,m_{k+1},\dots,m_{B_s}],
\]
and
\[
   E_k[w;m_1,\dots,m_{k-1},\nu,m_{k+1},\dots,m_{B_s}] =
   \bigg\{\prod_{i=1}^s \frac{q^{a_i
   m_{(1+B_{i-1})}}}{[m_{(1+B_{i-1})}]_q^{a_i}}\bigg\}
   \bigg/\prod_{\substack{j=1\\j\ne
   k}}^{B_s}\big([m_j]_q-q^{m_j-\nu}[\nu]_q\big).
\]
In the expression for $E_k$, we have placed the compound subscript
$1+B_{i-1}$ in parentheses to emphasize that the entire expression
$1+B_{i-1}$ occurs in the subscript of $m$.

\end{Lem}

We defer the proofs of Theorem~\ref{thm:qdiff} and
Lemma~\ref{lem:fmero} in order to proceed directly to the proof of
Theorem~\ref{thm:f=g}.

\noindent{\bf Proof of Theorem~\ref{thm:f=g}.}  We use induction
on the total degree of the word $\prod_{i=1}^s x^{a_i}y^{b_i}$.
The base case is clearly satisfied, since the word $xy$ is
self-dual. Now apply Theorem~\ref{thm:qdiff} to $f$ and $g$.
Subtracting the two equations gives
\begin{multline*}
   \sum_{\substack{\delta,\epsilon\in I^s\\ \delta_1<a_1, \epsilon_s<b_s}}
   (-\theta)^{\overline{\delta}\cdot\overline{\epsilon}}
   (1-q)^{\delta\cdot\epsilon}\bigg\{f\bigg[\prod_{i=1}^s
   x^{a_i-\delta_i}y^{b_i-\epsilon_i};\theta\bigg]-g\bigg[\prod_{i=1}^s
   x^{a_i-\delta_i}y^{b_i-\epsilon_i};\theta\bigg]\bigg\}\\
   = \sum_{\substack{\delta,\epsilon\in I^{s+1}\\
   \delta_{s+1}=\epsilon_1=0\\ \delta_1<a_1,\epsilon_{s+1}<b_s}}
   (-\theta\,')^{\overline{\delta}\cdot\overline{\epsilon}-1}(1-q)^{\delta\cdot\epsilon}
   \bigg\{f\bigg[\prod_{i=1}^s
   x^{a_i-\delta_i}y^{b_i-\epsilon_{i+1}};\theta\,'\bigg]
   -g\bigg[\prod_{i=1}^s
   x^{a_i-\delta_i}y^{b_i-\epsilon_{i+1}};\theta\,'\bigg]\bigg\}.
\end{multline*}
But the terms whose words have total degree less than
$\sum_{i=1}^s (a_i+b_i)$ are cancelled by the induction
hypothesis.  This leaves us with
\begin{multline*}
   (-\theta)^s\bigg\{f\bigg[\prod_{i=1}^s
   x^{a_i}y^{b_i};\theta\bigg]-g\bigg[\prod_{i=1}^s
   x^{a_i}y^{b_i};\theta\bigg]\bigg\}\\
   = (-\theta\,')^s\bigg\{f\bigg[\prod_{i=1}^s
   x^{a_i}y^{b_i};\theta\,'\bigg]-g\bigg[\prod_{i=1}^s
   x^{a_i}y^{b_i};\theta\,'\bigg]\bigg\}.
\end{multline*}
Thus, the function
\[
   H(\theta):=(-\theta)^s \bigg\{f\bigg[\prod_{i=1}^s
   x^{a_i}y^{b_i};\theta\bigg]-g\bigg[\prod_{i=1}^s
   x^{a_i}y^{b_i};\theta\bigg]\bigg\}
\]
satisfies the functional equation $H(\theta)=H(\theta\,')$, where
$\theta\,'=q\theta-1$.  But by Lemma~\ref{lem:fmero}, $H(\theta)$
is a meromorphic function of $\theta$ of the form
\[
   \theta^s\sum_{\nu=1}^\infty \frac{h_{\nu}}{[\nu]_q-\theta q^{\nu}},
\]
with at worst simple poles at $\theta=p_{\nu}:=q^{-\nu}[\nu]_q$
for positive integers $\nu$.  Note that $0=p_0<p_1<p_2<\cdots$ and
${p\,'}_{\!\!\!\!\nu}=qp_{\nu}-1=p_{\nu-1}$ for all $\nu\ge 1$.
The functional equation thus implies that if $H$ has a pole at
$p_{\nu}$, then $H$ must also have a pole at $p_{\nu-1}$.  Since
$H$ has no pole at $p_0$, it follows that each $h_{\nu}=0$.  Thus,
$H$ vanishes identically and the proof is complete. \eop

Let $1\ne w=\prod_{i=1}^s x^{a_i}y^{b_i}\in\mathfrak{h}^0$.
Henceforth, we assume that $|\theta|<1/q$.   To prove that $f$ and
$g$ satisfy the difference equation as stipulated by
Theorem~\ref{thm:qdiff}, first observe that from~\eqref{fgf},
\begin{align}
   f[w;\theta] &= \sum_{\nu=0}^\infty \theta^{\nu}
   \sum_{\substack{c_j\ge 0\\ \sum_{j=1}^n c_j=\nu}}\;
   \sum_{m_1>\cdots>m_n>0}\;\prod_{j=1}^n
   \frac{q^{(k_j+c_j-1)m_j}}{[m_j]_q^{k_j+c_j}}\notag\\
   &= \sum_{m_1>\cdots>m_n>0}\;\prod_{j=1}^n \sum_{c_j=0}^\infty
   \frac{q^{(k_j+c_j-1)m_j}}{[m_j]_q^{k_j+c_j}}\,\theta^{c_j}\notag\\
   &= \sum_{m_1>\cdots>m_n>0}\;\prod_{j=1}^n \frac{q^{(k_j-1)m_j}}
   {[m_j]_q^{k_j-1}\big([m_j]_q-\theta q^{m_j}\big)}\notag\\
   &= \sum_{m_1>\cdots>m_{B_s}>0}\; \prod_{i=1}^s
   \frac{q^{a_im_{(1+B_{i-1})}}}{[m_{(1+B_{i-1})}]_q^{a_i}}
   \prod_{j=1+B_{i-1}}^{B_i}\frac{1}{[m_j]_q-\theta q^{m_j}},
\label{fnest}
\end{align}
where $B_0:=0$ and $B_i:=\sum_{j=1}^i b_j$ for $1\le i\le s$ as in
the statement of Lemma~\ref{lem:fmero}.
\begin{Def}  If $d=(d_1,\dots,d_s)\in I^s$ is such that $d_s=0$
if $b_s=1$, let
\[
   f[w;d;\theta] := \sum_{m_1>\cdots>m_{B_s}>0}\;
   \prod_{i=1}^s \frac{q^{a_i(m_{(1+B_{i-1})}-d_i)}}
   {[m_{(1+B_{i-1})}-d_i]_q^{a_i}}
   \prod_{j=1+B_{i-1}}^{B_i}\frac{1}{[m_j]_q-\theta q^{m_j}}.
\]
\end{Def}
The extra requirement on $d_s$ ensures that no division by zero
occurs when $B_s=1$.  Note that we now have
$f[w;\theta]=f[w;\{0\}^s;\theta]$.  For the proof of
Theorem~\ref{thm:qdiff},  we require the following sequence of
lemmata.

\begin{Lem}\label{lem:14ia} If $(s>1$ or $b_1>1)$  and $a_1>1$, then
\begin{multline*}
   \sum_{\delta,\epsilon\in I}
   (-\theta)^{\overline{\delta}\overline{\epsilon}}(1-q)^{\delta\epsilon}
   f\big[x^{a_1-\delta}y^{b_1-\epsilon}\prod_{i=2}^s
   x^{a_i}y^{b_i};\{0\}^s;\theta\big]\\
   = \sum_{\delta\in I}(-\theta\,')^{\overline{\delta}}
   f\big[x^{a_1-\delta}y^{b_1}\prod_{i=2}^sx^{a_i}y^{b_i};1,\{0\}^{s-1};
   \theta\big].
\end{multline*}
\end{Lem}

\begin{Lem}\label{lem:14ib} If $(s>1$ or $b_1>1)$ and $a_1=1$ then
\[
   \sum_{\epsilon\in I} (-\theta)^{\overline{\epsilon}}
   f\big[xy^{b_1-\epsilon}\prod_{i=2}^s
   x^{a_i}y^{b_i};\{0\}^s;\theta\big]
   = (-\theta\,')f\big[xy^{b_1}\prod_{i=2}^s
   x^{a_i}y^{b_i};1,\{0\}^{s-1};\theta\big].
\]
\end{Lem}

\begin{Lem}\label{lem:14ii} If $1<j<s$ or $(j=s$ and $b_s>1)$ then
\begin{multline*}
   \sum_{\delta,\epsilon\in I}
   (-\theta)^{\overline{\delta}\overline{\epsilon}} (1-q)^{\delta\epsilon}
   f\bigg[\big(\prod_{i=1}^{j-1}x^{a_i}y^{b_i}\big)x^{a_j-\delta}
   y^{b_j-\epsilon}\prod_{i=j+1}^s x^{a_i}y^{b_i};
   \{1\}^{j-1},\{0\}^{s-j+1};\theta\bigg]\\
   = \sum_{\delta,\epsilon\in I}
   (-\theta\,')^{\overline{\delta}\overline{\epsilon}} (1-q)^{\delta\epsilon}
   f\bigg[\big(\prod_{i=1}^{j-2}x^{a_i}y^{b_i}\big)x^{a_{j-1}}
   y^{b_{j-1}-\epsilon}x^{a_j-\delta}y^{b_j}
   \prod_{i=j+1}^s x^{a_i}y^{b_i};
   \{1\}^j,\{0\}^{s-j};\theta\bigg].
\end{multline*}
\end{Lem}

\begin{Lem}\label{lem:14iiia} If $b_s>1$, then
\[
   f\bigg[\prod_{i=1}^s x^{a_i}y^{b_i};\{1\}^s;\theta\bigg]
   = \sum_{\epsilon\in I} (-\theta\,')^{-\epsilon}
   f\bigg[\big(\prod_{i=1}^{s-1}x^{a_i}y^{b_i}\big)x^{a_s}y^{b_s-\epsilon};
   \{0\}^s;\theta\,'\bigg].
\]
\end{Lem}

\begin{Lem}\label{lem:14iiib} If $s>1$, then
\begin{multline*}
   \sum_{\delta\in I} (-\theta)^{\overline{\delta}} f\bigg[
   \big(\prod_{i=1}^{s-1}x^{a_i}y^{b_i}\big)x^{a_s-\delta}y;\{1\}^{s-1};
   0;\theta\bigg]\\
   =\sum_{\delta,\epsilon\in I}(-\theta\,')^{\overline{\delta}\overline{\epsilon}}
   (1-q)^{\delta\epsilon}f\bigg[\big(\prod_{i=1}^{s-2}x^{a_i}y^{b_i}\big)
   x^{a_{s-1}}y^{b_{s-1}-\epsilon}x^{a_s-\delta}y;\{0\}^s;\theta\,'\bigg].
\end{multline*}
\end{Lem}

\begin{Lem}\label{lem:14iiic} If $a>1$ then
\[
   \sum_{\delta\in I} (-\theta)^{\overline{\delta}}
   f\big[x^{a-\delta}y;\theta\big]
   = \sum_{\delta\in I} (-\theta\,')^{\overline{\delta}}
   f\big[x^{a-\delta}y;\theta\,'\big].
\]
\end{Lem}

For completeness, we also record the following result, although it
is not needed for the proof of Theorem~\ref{thm:qdiff}.
\begin{Lem}\label{lem:14iiid}
$
  \theta f[xy;\theta] +(1-q) = \theta\,'f[xy;\theta\,']
   -1/\theta\,'.
$
\end{Lem}

We shall prove Lemmas~\ref{lem:fmero}--\ref{lem:14iiid} in
Subsection~\ref{sect:lemmas} below.  Assuming their validity for
now, we proceed with the proof of Theorem~\ref{thm:qdiff}.

\noindent{\bf Proof of Theorem~\ref{thm:qdiff}.} Let $L$ denote
the left hand side.  First, consider the case when $a_1>1$ and
$b_s>1$. Then
\[
   L = \sum_{\delta,\epsilon\in I^s} (-\theta)^{\overline{\delta}
   \cdot\overline{\epsilon}} (1-q)^{\delta\cdot\epsilon}
   f\bigg[\prod_{i=1}^s
   x^{a_i-\delta_i}y^{b_i-\epsilon_i};\theta\bigg].
\]
In the sum over ordered $s$-tuples $\delta$ and $\epsilon$, rename
$\delta=(\delta_2,\dots,\delta_s)$ and
$\epsilon=(\epsilon_2,\dots,\epsilon_s)$ so that
\begin{align*}
   L &= \sum_{\delta,\epsilon\in I^{s-1}}
   (-\theta)^{\overline{\delta}\cdot\overline{\epsilon}}
   (1-q)^{\delta\cdot\epsilon}
   \sum_{\delta_1,\epsilon_1\in I}(-\theta)^{\overline{\delta}_1
   \overline{\epsilon}_1}(1-q)^{\delta_1\epsilon_1}
   f\bigg[\prod_{i=1}^s
   x^{a_i-\delta_i}y^{b_i-\epsilon_i};\{0\}^s;\theta\bigg]\\
   &= \sum_{\delta,\epsilon\in I^{s-1}}
    (-\theta)^{\overline{\delta}\cdot\overline{\epsilon}}
   (1-q)^{\delta\cdot\epsilon}
   \sum_{\delta_1\in I}(-\theta\,')^{\overline{\delta}_1}
    f\bigg[x^{a_1-\delta_1}y^{b_1}\prod_{i=2}^s
   x^{a_i-\delta_i}y^{b_i-\epsilon_i};1,\{0\}^{s-1};\theta\bigg],
\end{align*}
by Lemma~\ref{lem:14ia}.  If $s>1$, we again rename
$\delta=(\delta_3,\dots,\delta_s)$ and
$\epsilon=(\epsilon_3,\dots,\epsilon_s)$ and write
\begin{multline*}
   L = \sum_{\delta_1\in I}(-\theta\,')^{\overline{\delta}_1}
   \sum_{\delta,\epsilon\in I^{s-2}}(-\theta)^{\overline{\delta}
   \cdot\overline{\epsilon}}(1-q)^{\delta\cdot\epsilon}\\
   \times \sum_{\delta_2,\epsilon_2\in I}
   (-\theta)^{\overline{\delta}_2\overline{\epsilon}_2}
   (1-q)^{\delta_2\epsilon_2}
   f\bigg[x^{a_1-\delta_1}y^{b_1}\prod_{i=2}^s x^{a_i-\delta_i}
   y^{b_i-\epsilon_i};1,\{0\}^{s-1};\theta\bigg].
\end{multline*}
We now apply Lemma~\ref{lem:14ii}, first with $j=2$, and again
with $j=3$, and so on up to $j=s$.  The result is that
\begin{equation}\label{no-app}
   L=\sum_{\substack{\delta=(\delta_1,\dots,\delta_s)\in I^s\\
   \epsilon=(\epsilon_1,\dots,\epsilon_s)\in I^s\\
   \epsilon_1=0}}
   (-\theta\,')^{\overline{\delta}\cdot\overline{\epsilon}}
   (1-q)^{\delta\cdot\epsilon}f\bigg[\big(\prod_{i=1}^{s-1}
   x^{a_i-\delta_i}y^{b_i-\epsilon_{i+1}}\big)x^{a_s-\delta_s}y^{b_s};
   \{1\}^s;\theta\bigg].
\end{equation}
On the other hand, if $s=1$, we have~\eqref{no-app} with no
application of Lemma~\ref{lem:14ii}.  In any case, applying
Lemma~\ref{lem:14iiia} to~\eqref{no-app} yields
\[
   L=\sum_{\substack{\delta,\epsilon\in I^s\\ \epsilon_1=0}}
   (-\theta\,')^{\overline{\delta}\cdot\overline{\epsilon}}
   (1-q)^{\delta\cdot\epsilon}
   \sum_{\epsilon_{s+1}\in I}(-\theta\,')^{\overline{\epsilon}_{s+1}-1}
   f\bigg[\prod_{i=1}^s
   x^{a_i-\delta_i}y^{b_i-\epsilon_{i+1}};\{0\}^s;\theta\,'\bigg].
\]
If we now extend $\delta$ and $\epsilon$ by adjoining an extra
component to each, viz.\ $\delta_{s+1}=0$ and $\epsilon_{s+1}\in
I$ respectively, we find that
\[
   L= \sum_{\substack{\delta,\epsilon\in I^{s+1}\\
   \delta_{s+1}=\epsilon_1=0}} (-\theta\,')^{\overline{\delta}\cdot
   \overline{\epsilon}}(1-q)^{\delta\cdot\epsilon}
   f\bigg[\prod_{i=1}^s
   x^{a_i-\delta_i}y^{b_i-\epsilon_{i+1}};\{0\}^s;\theta\,'\bigg],
\]
as required.

The proof in the case $a_1=1$, $b_s>1$ is similar.  The main
difference is that $\delta_1=0$ and we begin by applying
Lemma~\ref{lem:14ib} instead of Lemma~\ref{lem:14ia}.  For
purposes of brevity, we suppress the details.

It is convenient to split the case $a_1>1$, $b_s=1$ into the two
subcases $s>1$ and $s=1$, since in the former we end by applying
Lemma~\ref{lem:14iiib}, while in the latter we instead use
Lemma~\ref{lem:14iiic}.  Suppose first that $s>1$.  We have
\begin{align*}
   L&=\sum_{\substack{\delta,\epsilon\in I^s\\ \epsilon_s=0}}
   (-\theta)^{\overline{\delta}\cdot\overline{\epsilon}}(1-q)^{\delta\cdot\epsilon}
   f\bigg[\prod_{i=1}^s x^{a_i-\delta_i}
   y^{b_i-\epsilon_i};\{0\}^s;\theta\bigg]\\
   &= \sum_{\substack{\delta,\epsilon\in I^{s-1}\\
   \epsilon_s=0}}
   (-\theta)^{\overline{\delta}\cdot\overline{\epsilon}}
   (1-q)^{\delta\cdot\epsilon}\sum_{\delta_1,\epsilon_1\in I}
   f\bigg[x^{a_1-\delta_1}y^{b_1-\epsilon_1}\prod_{i=2}^s
   x^{a_i-\delta_i}y^{b_i-\epsilon_i};\{0\}^s;\theta\bigg].
\end{align*}
Now apply Lemma~\ref{lem:14ia}, and then Lemma~\ref{lem:14ii}
successively, with $j=2,3,\dots,s-1$.  The result is
\[
   L= \sum_{\substack{\eta,\nu\in I^{s-1}\\
   \nu_1=0}}(-\theta\,')^{\overline{\eta}\cdot\overline{\nu}}
   (1-q)^{\eta\cdot\nu}\sum_{\substack{\delta_s\in I\\ \nu_s=0}}(-\theta)^{\overline{\delta}_s}
   f\bigg[\big(\prod_{i=1}^{s-1}x^{a_i-\eta_i}y^{b_i-\nu_{i+1}}\big)
   x^{a_s-\delta_s}y;\{1\}^{s-1},0;\theta\bigg].
\]
Lemma~\ref{lem:14iiib} now gives
\begin{align*}
   L&=\sum_{\substack{\eta,\nu\in I^s\\ \nu_1=0}}
   (-\theta\,')^{\overline{\eta}\cdot\overline{\nu}}(1-q)^{\eta\cdot\nu}
   f\bigg[\prod_{i=1}^{s-1}
   x^{a_i-\eta_i}y^{b_i-\nu_{i+1}}x^{a_s-\eta_s}y;\{0\}^s;\theta\,'\bigg]\\
   &= \sum_{\substack{\eta,\nu\in I^{s+1}\\
   \eta_{s+1}=\nu_1=\nu_{s+1}=0}}
   (-\theta\,')^{\overline{\eta}\cdot\overline{\nu}-1}(1-q)^{\eta\cdot\nu}
   f\bigg[\prod_{i=1}^s
   x^{a_i-\eta_i}y^{b_i-\nu_{i+1}};\theta\,'\bigg],
\end{align*}
as required.  On the other hand, if $s=1$ note that in this case
Theorem~\ref{thm:qdiff} is just a restatement of
Lemma~\ref{lem:14iiic}.

The final case, with $a_1=b_s=1$ and $s>1$, is proved in much the
same way as the other cases with $s>1$.  Observe that now
$\delta_1=\epsilon_s=0$ in the sum on the left, and
$\delta_1=\epsilon_{s+1}=0$ on the right.  The result is
established by applying Lemma~\ref{lem:14ib}, then
Lemma~\ref{lem:14ii} successively as necessary for
$j=2,3,\dots,s-1$, and finally Lemma~\ref{lem:14iiib}.

Thus, $f$ satisfies the difference equation as claimed.  This and
the fact that $g[w;\theta]=f[\tau(w);\theta]$ readily implies that
$g$ satisfies the same difference equation.

\eop

\subsection{Proofs of Lemmas~\ref{lem:fmero}--\ref{lem:14iiid}}
\label{sect:lemmas}

We begin with the proof of Lemma~\ref{lem:fmero}. From the
penultimate step in~\eqref{fnest}, noting that $n=B_s$, we have
\begin{align*}
   f[w;\theta] &= \sum_{m_1>\cdots>m_{B_s}>0}\;
   \prod_{j=1}^{B_s} \frac{q^{(k_j-1)m_j}}{[m_j]_q^{k_j-1}
   \big([m_j]_q-\theta q^{m_j}\big)}
   = \sum_{m_1>\cdots>m_{B_s}>0}\;
   \sum_{k=1}^{B_s} \frac{E_k[w;m_1,\dots,m_{B_s}]}
   {[m_k]_q-\theta q^{m_k}},
\end{align*}
where the partial fraction decomposition
\[
   \sum_{h=1}^{B_s} \frac{E_h[w;m_1,\dots,m_{B_s}]}
   {[m_h]_q-\theta q^{m_h}}
   = \prod_{j=1}^{B_s}\frac{q^{(k_j-1)m_j}}{[m_j]_q^{k_j-1}
   \big([m_j]_q-\theta q^{m_j}\big)}
\]
implies that
\[
   \sum_{h=1}^{B_s}E_h[w;m_1,\dots,m_{B_s}]
   \prod_{\substack{j=1\\ j\ne h}}^{B_s}\big([m_j]_q-\theta
   q^{m_j}\big) =
   \prod_{j=1}^{B_s} \frac{q^{(k_j-1)m_j}}{[m_j]_q^{k_j-1}}.
\]
Letting $\theta\to q^{-m_k}[m_k]_q$ now gives that
\[
  E_k[w;m_1,\dots,m_{B_s}] =
  \bigg\{\prod_{j=1}^{B_s} \frac{q^{(k_j-1)m_j}}{[m_j]_q^{k_j-1}}
  \bigg\}\bigg/\prod_{\substack{j=1\\ j\ne k}}^{B_s}
  \big([m_j]_q-q^{m_j-m_k}[m_k]_q\big).
\]
The general formula for
$E_k[m_1,\dots,m_{k-1},\nu,m_{k+1},\dots,m_{B_s}]$ now follows
immediately on replacing $m_k$ by $\nu$ and noting that
$k_j=a_i+1$ precisely when $j=1+B_{i-1}$; otherwise $k_j=1$. The
lemma itself now follows on interchanging order of summation. \eop

Proofs of several of the remaining lemmata make use of the partial
fraction identity
\begin{multline}\label{parfrac}
   \frac{\theta q^{2m}}{[m]_q^a\big([m]_q-\theta q^m\big)}
   -\frac{q^m}{[m]_q^{a-1}\big([m]_q-\theta q^m\big)}\\
   = \frac{\theta\,' q^{2m-a}}{[m-1]_q^a\big([m]_q-\theta
   q^m\big)}-\frac{q^{m-a+1}}{[m-1]_q^{a-1}\big([m]_q-\theta
   q^m\big)} + \frac{q^{m-a+1}}{[m-1]_q^a}-\frac{q^m}{[m]_q^a},
\end{multline}
valid for $a>0$ and $m>1$.

\noindent{\bf Proof of Lemma~\ref{lem:14ia}.} Let
\begin{equation}\label{B}
   B:= \prod_{j=2}^n
   \frac{q^{(k_j-1)m_j}}{[m_j]_q^{k_j-1}\big([m_j]_q-\theta
   q^{m_j}\big)}.
\end{equation}
Then by~\eqref{parfrac},
\begin{align*}
   &\theta f\bigg[\prod_{i=1}^s x^{a_i}
   y^{b_i};\{0\}^s;\theta\bigg] -
   f\bigg[x^{a_1-1}y^{b_1}\prod_{i=2}^s
   x^{a_i}y^{b_i};\{0\}^s;\theta\bigg]\\
   &=\sum_{m_1>\cdots>m_n>0}\;
   \bigg\{\frac{\theta q^{2m_1}}{[m_1]_q^{a_1}
   \big([m_1]_q-\theta q^{m_1}\big)}
   -\frac{q^m_1}{[m_1]_q^{a_1-1}\big([m_1]_q-\theta
   q^{m_1}\big)}\bigg\}
   q^{(a_1-2)m_1}B\\
   &=\sum_{m_1>\cdots>m_n>0}\;\bigg\{\frac{\theta\,' q^{2m_1-a_1}}
   {[m_1-1]_q^{a_1}\big([m_1]_q-\theta
   q^{m_1}\big)}-\frac{q^{m_1-a_1+1}}{[m_1-1]_q^{a_1-1}\big([m_1]_q-\theta
   q^{m_1}\big)}\\
   &\qquad\qquad\qquad\qquad\qquad\qquad\qquad\qquad\qquad\qquad
    + \frac{q^{m_1-a_1+1}}{[m_1-1]_q^{a_1}}
   -\frac{q^{m_1}}{[m_1]_q^{a_1}}\bigg\}q^{(a_1-2)m_1}B\\
   &=\sum_{m_1>\cdots>m_n>0}\;\bigg\{
   \frac{\theta\,'q^{a_1(m_1-1)}}{[m_1-1]_q^{a_1}\big([m_1]_q-
   \theta
   q^{m_1}\big)}-\frac{q^{(a_1-1)(m_1-1)}}{[m_1-1]_q^{a_1-1}
   \big([m_1]_q-\theta q^{m_1}\big)}\bigg\}B\\
   &\qquad+\sum_{m_2>\dots>m_n>0}\; B\sum_{m_1=m_2+1}^\infty
   \bigg\{\frac{q^{(a_1-1)(m_1-1)}}{[m_1-1]_q^{a_1}}-
   \frac{q^{(a_1-1)m_1}}{[m_1]_q^{a_1}}\bigg\}\\
   &= \theta\,'f\bigg[\prod_{i=1}^s
   x^{a_i}y^{b_i};1;\{0\}^{s-1}\;\theta\bigg]
   - f\bigg[x^{a_1-1}y^{b_1}\prod_{i=2}^s x^{a_i}y^{b_i};1;
   \{0\}^{s-1};\theta\bigg]\\
   &\qquad + \sum_{m_2>\cdots>m_n>0}
   \frac{Bq^{(a_1-1)m_2}}{[m_2]_q^{a_1}}.
\end{align*}
But
\begin{align*}
  &\sum_{m_2>\cdots>m_n>0}
   \frac{Bq^{(a_1-1)m_2}}{[m_2]_q^{a_1}}
   = \sum_{m_2>\cdots>m_n>0}\bigg\{
   \frac{q^{a_1m_2}}{[m_2]_q^{a_1}}+
   (1-q)\frac{q^{(a_1-1)m_2}}{[m_2]_q^{a_1-1}}\bigg\}B\\
   &=f\bigg[x^{a_1}y^{b_1-1}\prod_{i=2}^s x^{a_i}y^{b_i};
   \{0\}^{s};\theta\bigg]
   +(1-q)f\bigg[x^{a_1-1}y^{b_1-1}\prod_{i=2}^s x^{a_i}y^{b_i};
   \{0\}^{s};\theta\bigg],
\end{align*}
and the result follows. \eop

\noindent{\bf Proof of Lemma~\ref{lem:14ib}.} Again, let $B$ be
given by~\eqref{B}.  In this case~\eqref{parfrac} gives
\begin{align*}
   &\theta f\bigg[xy^{b_1}\prod_{i=2}^s
   x^{a_i}y^{b_i};\{0\}^s;\theta\bigg]
   =\sum_{m_1>\cdots>m_n>0}\frac{\theta q^{2m_1}}{[m_1]_q\big(
   [m_1]_q-\theta q^{m_1}\big)}\cdot q^{-m_1}B\\
   &=\sum_{m_1>\cdots>m_n>0}\bigg\{\frac{\theta\,' q^{2m_1-1}}
   {[m_1-1]_q\big([m_1]_q-\theta q^{m_1}\big)}
   +\frac{q^{m_1}}{[m_1-1]_q}-\frac{q^{m_1}}{[m_1]_q}\bigg\}
   q^{-m_1}B\\
   &=\sum_{m_1>\cdots>m_n>0}\frac{\theta\,' q^{m_1-1}B}{[m_1-1]_q
   \big([m_1]-\theta q^{m_1}\big)}+
   \sum_{m_2>\cdots>m_n>0} B\sum_{m_1=m_2+1}^\infty
   \bigg(\frac{1}{[m_1-1]_q}-\frac{1}{[m_1]_q}\bigg)\\
   &= \theta\,'f\bigg[xy^{b_1}\prod_{i=2}^sx^{a_i}y^{b_i};1,
   \{0\}^{s-1};\theta\bigg]
   +\sum_{m_2>\cdots>m_n>0}\bigg(\frac{1}{[m_2]_q}+q-1\bigg)B.
\end{align*}
In light of $q-1+1/[m_2]_q=q^{m_2}/[m_2]_q$, it follows that
\begin{multline*}
   \theta f\bigg[xy^{b_1}\prod_{i=2}^s
   x^{a_i}y^{b_i};\{0\}^s;\theta\bigg]-
   \theta\,'f\bigg[xy^{b_1}\prod_{i=2}^sx^{a_i}y^{b_i};1,
   \{0\}^{s-1};\theta\bigg]\\
   = \sum_{m_2>\cdots>m_n>0}\frac{Bq^{m_2}}{[m_2]_q}
   =f\bigg[xy^{b_1-1}\prod_{i=2}^sx^{a_i}y^{b_i};\{0\}^s;\theta\bigg],
\end{multline*}
as claimed. \eop

\noindent{\bf Proof of Lemma~\ref{lem:14ii}.}  Let
$m=m_{(1+B_{j-1})}$.  Define the quantities $A$ and $B$ by
\[
   A = \prod_{i=1}^{j-1} \frac{q^{a_i(m_{(1+B_{i-1})}-d_i)}}
   {[m_{(1+B_{i-1})}-d_i]_q^{a_i}}\prod_{h=1+B_{i-1}}^{B_i}
   \frac{1}{[m_h]_q-\theta q^{m_h}},
\]
and
\[
   \frac{q^{a_jm}B}{[m]_q^{a_j}\big([m]_q-\theta q^m\big)}
   =\prod_{i=j}^s
   \frac{q^{a_im_{(1+B_{i-1})}}}{[m_{(1+B_{i-1})}]_q^{a_i}}
   \prod_{h=1+B_{i-1}}^{B_i}\frac{1}{[m_h]_q-\theta q^{m_h}}.
\]
Then~\eqref{parfrac} gives
\begin{align*}
   &\theta f\bigg[\prod_{i=1}^s
   x^{a_i}y^{b_i};\{1\}^{j-1},\{0\}^{s-j+1};\theta\bigg]\\
   &\qquad -f\bigg[\big(\prod_{i=1}^{j-1}x^{a_i}y^{b_i}\big)x^{a_j-1}
   y^{b_j}\prod_{i=j+1}^s x^{a_i}y^{b_i};
   \{1\}^{j-1},\{0\}^{s-j+1};\theta\bigg]\\
   &=\sum_{m_1>\cdots>m_n>0}
   A\bigg\{\frac{\theta q^{2m}}{[m]_q^{a_j}\big([m]_q-\theta q^m\big)}-
   \frac{q^m}{[m]_q^{a_j-1}\big([m]_q-\theta q^m\big)}\bigg\}
   q^{(a_j-2)m}B\\
   &= \sum_{m_1>\cdots>m_n>0}A\bigg\{\frac{\theta\,' q^{2m-a_j}}
   {[m-1]_q^{a_j}\big([m]_q^{a_j}-\theta q^{m}\big)}
   -\frac{q^{m-a_j+1}}{[m-1]_q^{a_j-1}\big([m]_q-\theta
   q^m\big)}\\
   &\qquad\qquad\qquad\qquad\qquad\qquad
   +\frac{q^{m-a_j+1}}{[m-1]_q^{a_j}}-\frac{q^m}{[m]_q^{a_j}}\bigg\}
   q^{(a_j-2)m}B\\
   &=\sum_{m_1>\cdots>m_n>0} A\bigg\{\frac{\theta\,' q^{a_j(m-1)}}
   {[m-1]_q^{a_j}\big([m]_q-\theta q^m\big)}-
   \frac{q^{(a_j-1)(m-1)}}{[m-1]_q^{a_j-1}\big([m]_q-\theta
   q^m\big)}\bigg\}B\\
   &\qquad\qquad +\sum_{m_1>\cdots>m_n>0} A\bigg\{
   \frac{q^{(a_j-1)(m-1)}}{[m-1]_q^{a_j}}-\frac{q^{(a_j-1)m}}{[m]_q^{a_j}}
   \bigg\}B\\
   &= \theta\,' f\bigg[\prod_{i=1}^s
   x^{a_i}y^{b_i};\{1\}^j,\{0\}^{s-j};\theta\bigg]
   -f\bigg[\big(\prod_{i=1}^{j-1}x^{a_i}y^{b_i}\big)x^{a_j-1}y^{b_j}
   \prod_{i=j+1}^s
   x^{a_i}y^{b_i};\{1\}^j,\{0\}^{s-j};\theta\bigg]\\
   &\qquad + \sum_{\substack{m_1>\cdots>m_{B_{j-1}}\\
   -1+m_{B_{j-1}}>m_{(2+B_{j-1})}>\cdots>m_{B_s}>0}}AB
   \sum_{m=1+m_{(2+B_{j-1})}}^{-1+m_{B_{j-1}}}
   \bigg(\frac{q^{(a_j-1)(m-1)}}{[m-1]_q^{a_j}}-
   \frac{q^{(a_j-1)m}}{[m]_q^{a_j}}\bigg).
\end{align*}
It follows that
\begin{align*}
   &\theta f\bigg[\prod_{i=1}^s
   x^{a_i}y^{b_i};\{1\}^{j-1},\{0\}^{s-j+1};\theta\bigg]\\
   &-f\bigg[\big(\prod_{i=1}^{j-1}x^{a_i}y^{b_i}\big)x^{a_j-1}
   y^{b_j}\prod_{i=j+1}^s x^{a_i}y^{b_i};
   \{1\}^{j-1},\{0\}^{s-j+1};\theta\bigg]\\
   &-\theta\,' f\bigg[\prod_{i=1}^s
   x^{a_i}y^{b_i};\{1\}^j,\{0\}^{s-j};\theta\bigg]
   +f\bigg[\big(\prod_{i=1}^{j-1}x^{a_i}y^{b_i}\big)x^{a_j-1}y^{b_j}
   \prod_{i=j+1}^s
   x^{a_i}y^{b_i};\{1\}^j,\{0\}^{s-j};\theta\bigg]\\
   &= \sum_{\substack{m_1>\cdots>m_{B_{j-1}}\\
   -1+m_{B_{j-1}}>m_{(2+B_{j-1})}>\cdots>m_{B_s}>0}}A
   \bigg\{\frac{q^{(a_j-1)m_{(2+B_{j-1})}}}{[m_{(2+B_{j-1})}]_q^{a_j}}
   -\frac{q^{(a_j-1)(-1+m_{B_{j-1}})}}{[-1+m_{B_{j-1}}]_q^{a_j}}\bigg\}B\\
   &=\sum_{m_1>\cdots>m_{B_{j-1}}>m_{(2+B_{j-1})}>\cdots>m_{B_s}>0}
   A\bigg\{\frac{q^{a_jm_{(2+B_{j-1})}}}{[m_{(2+B_{j-1})}]_q^{a_j}}+
   (1-q)\frac{q^{(a_j-1)m_{(2+B_{j-1})}}}{[m_{(2+B_{j-1})}]_q^{a_j-1}}\\
   &\qquad\qquad\qquad\qquad\qquad\qquad\qquad
    -\frac{q^{a_j(-1+m_{B_{j-1}})}}{[-1+m_{B_{j-1}}]_q^{a_j}}
   -(1-q)\frac{q^{(a_j-1)(-1+m_{B_{j-1}})}}{[-1+m_{B_{j-1}}]_q^{a_j-1}}
  \bigg\}B\\
  &=f\bigg[\big(\prod_{i=1}^{j-1}x^{a_i}y^{b_i}\big)x^{a_j}y^{b_j-1}
  \prod_{i=j+1}^s
  x^{a_i}y^{b_i};\{1\}^{j-1},\{0\}^{s-j+1};\theta\bigg]\\
  &\quad+(1-q)f\bigg[\big(\prod_{i=1}^{j-1}x^{a_i}y^{b_i}\big)x^{a_j-1}
  y^{b_j-1}\prod_{i=j+1}^s
  x^{a_i}y^{b_i};\{1\}^{j-1},\{0\}^{s-j+1};\theta\bigg]\\
  &\quad-f\bigg[\big(\prod_{i=1}^{j-2} x^{a_i}y^{b_i}\big)x^{a_{j-1}}
  y^{b_{j-1}-1}\prod_{i=j}^s x^{a_i}y^{b_i};\{1\}^j,
  \{0\}^{s-j};\theta\bigg]\\
  &\quad-(1-q)f\bigg[\big(\prod_{i=1}^{j-2}x^{a_j}y^{b_j}\big)
  x^{a_{j-1}}y^{b_{j-1}-1}x^{a_j-1}y^{b_j}\prod_{i=j+1}^s
  x^{a_i}y^{b_i};\{1\}^j,\{0\}^{s-j};\theta\bigg],
\end{align*}
as required. \eop

\noindent{\bf Proof of Lemma~\ref{lem:14iiia}.} Here $b_s>1$, and
thus if we shift summation indices $m_i\mapsto 1+m_i$, then
\begin{align*}
   &f\bigg[\prod_{i=1}^s x^{a_i}y^{b_i};\{1\}^s;\theta\bigg]
   =\sum_{m_1>\cdots>m_{B_s}>0}\;\prod_{i=1}^s
   \frac{q^{a_i(m_{(1+B_{i-1})}-1)}}{[m_{(1+B_{i-1})}-1]_q^{a_i}}
   \prod_{j=1+B_{i-1}}^{B_i}\frac{1}{[m_j]_q-\theta q^{m_j}}\\
   &=\sum_{m_1>\cdots>m_{B_s}\ge 0}\; \prod_{i=1}^s
   \frac{q^{a_im_{(1+B_{i-1})}}}{[m_{(1+B_{i-1})}]_q^{a_i}}
   \prod_{j=1+B_{i-1}}^{B_i}\frac{1}{[m_j+1]_q-\theta
   q^{m_j+1}}\\
   &= \sum_{m_1>\cdots>m_{B_s}\ge 0}\;\prod_{i=1}^s
   \frac{q^{a_im_{(1+B_{i-1})}}}{[m_{(1+B_{i-1})}]_q^{a_i}}
   \prod_{j=1+B_{i-1}}^{B_i}\frac{1}{[m_j]_q-\theta\,'q^{m_j}}\\
   &= \bigg(\sum_{m_1>\cdots>m_{B_s}> 0}+\sum_{m_1>\cdots>m_{(B_s-1)}>0}
   \bigg)\;\prod_{i=1}^s
   \frac{q^{a_im_{(1+B_{i-1})}}}{[m_{(1+B_{i-1})}]_q^{a_i}}
   \prod_{j=1+B_{i-1}}^{B_i}\frac{1}{[m_j]_q-\theta\,'q^{m_j}}\\
   &= f\bigg[\prod_{i=1}^s x^{a_i}y^{b_i};\{0\}^s;\theta\,'\bigg]
  -\bigg(\frac{1}{\theta\,'}\bigg)
  f\bigg[\big(\prod_{i=1}^{s-1}x^{a_i}y^{b_i}\big)
   x^{a_s}y^{b_s-1};\{0\}^s;\theta\,'\bigg].
\end{align*}
\eop

\noindent{\bf Proof of Lemma~\ref{lem:14iiib}.} In this case,
$B_s=1+B_{s-1}$ and we have
\begin{align*}
   &f\bigg[\big(\prod_{i=1}^{s-1}x^{a_i}y^{b_i}\big)x^{a_s-1}y;
   \{1\}^{s-1},0;\theta\bigg]
   -\theta
   f\bigg[\big(\prod_{i=1}^{s-1}x^{a_i}y^{b_i}\big)x^{a_s}y;
   \{1\}^{s-1},0;\theta\bigg]\\
   &= \sum_{m_1>\cdots>m_{B_s}>0}\;
   \bigg\{\prod_{i=1}^{s-1}\frac{q^{a_i(m_{(1+B_{i-1})}-1)}}
   {[m_{(1+B_{i-1})}-1]_q^{a_i}}\prod_{j=1+B_{i-1}}^{B_i}
   \frac{1}{[m_j]_q-\theta q^{m_j}}\bigg\}\\
   &\qquad\qquad\qquad\qquad\qquad \times
   \bigg\{\frac{q^{(a_s-1)m_{B_s}}}{[m_{B_s}]_q^{a_s-1}}-\frac{\theta
   q^{a_sm_{B_s}}}{[m_{B_s}]_q^{a_s}}
   \bigg\}\frac{1}{[m_{B_s}]_q-\theta
   q^{m_{B_s}}}\\
   &= \sum_{m_1>\cdots>m_{B_s}>0}\bigg\{\prod_{i=1}^{s-1}\frac{q^{a_i(m_{(1+B_{i-1})}-1)}}
   {[m_{(1+B_{i-1})}-1]_q^{a_i}}\prod_{j=1+B_{i-1}}^{B_i}
   \frac{1}{[m_j]_q-\theta
   q^{m_j}}\bigg\}\frac{q^{(a_s-1)m_{B_s}}}{[m_{B_s}]_q^{a_s}}.
\end{align*}
Now shift the summation indices $m_i\mapsto m_i+1$ and use
$[m+1]_q-\theta q^{m+1}=[m]_q-\theta\,' q^m$ to obtain
\begin{align*}
   &f\bigg[\big(\prod_{i=1}^{s-1}x^{a_i}y^{b_i}\big)x^{a_s-1}y;
   \{1\}^{s-1},0;\theta\bigg]
   -\theta
   f\bigg[\big(\prod_{i=1}^{s-1}x^{a_i}y^{b_i}\big)x^{a_s}y;
   \{1\}^{s-1},0;\theta\bigg]\\
   &= \sum_{m_1>\cdots>m_{B_s}\ge 0}\bigg\{
   \prod_{i=1}^{s-1}\frac{q^{a_im_{(1+B_{i-1})}}}
   {[m_{(1+B_{i-1})}]_q^{a_i}}\prod_{j=1+B_{i-1}}^{B_i}
   \frac{1}{[m_j]_q-\theta\,'
   q^{m_j}}\bigg\}\frac{q^{(a_s-1)m_{1+B_s}}}{[1+m_{B_s}]_q^{a_s}}.\\
\end{align*}
Now replace $m_{B_s}$ by $m_{B_s}-1$.  Then
\begin{align*}
    &f\bigg[\big(\prod_{i=1}^{s-1}x^{a_i}y^{b_i}\big)x^{a_s-1}y;
   \{1\}^{s-1},0;\theta\bigg]
   -\theta
   f\bigg[\big(\prod_{i=1}^{s-1}x^{a_i}y^{b_i}\big)x^{a_s}y;
   \{1\}^{s-1},0;\theta\bigg]\\
   &=\sum_{m_1>\cdots>m_{B_s-1}\ge m_{B_s}>0}\bigg\{
   \prod_{i=1}^{s-1}\frac{q^{a_im_{(1+B_{i-1})}}}
   {[m_{(1+B_{i-1})}]_q^{a_i}}\prod_{j=1+B_{i-1}}^{B_i}
   \frac{1}{[m_j]_q-\theta\,'
   q^{m_j}}\bigg\}\frac{q^{(a_s-1)m_{B_s}}}{[m_{B_s}]_q^{a_s}}\\
   &=\sum_{m_1>\cdots>m_{B_s-1}> m_{B_s}>0}\bigg\{
   \prod_{i=1}^{s-1}\frac{q^{a_im_{(1+B_{i-1})}}}
   {[m_{(1+B_{i-1})}]_q^{a_i}}\prod_{j=1+B_{i-1}}^{B_i}
   \frac{1}{[m_j]_q-\theta\,'
   q^{m_j}}\bigg\}\\
   &\qquad\qquad\qquad\qquad\qquad\qquad\qquad\qquad\qquad
   \times\frac{q^{(a_s-1)m_{B_s}}}{[m_{B_s}]_q^{a_s}}
   \cdot\frac{[m_{B_s}]_q-\theta\,'q^{m_{B_s}}}
   {[m_{B_s}]_q-\theta\,'q^{m_{B_s}}}\\
   &\quad+\sum_{m_1>\cdots>m_{B_s-1}>0}\bigg\{
   \prod_{i=1}^{s-1}\frac{q^{a_im_{(1+B_{i-1})}}}
   {[m_{(1+B_{i-1})}]_q^{a_i}}\prod_{j=1+B_{i-1}}^{B_i}
   \frac{1}{[m_j]_q-\theta\,'
   q^{m_j}}\bigg\}\frac{q^{(a_s-1)m_{({B_s}-1)}}}{[m_{(B_s-1)}]_q^{a_s}}\\
   &=\sum_{m_1>\cdots>m_{B_s}>0}\bigg\{
   \prod_{i=1}^{s-1}\frac{q^{a_im_{(1+B_{i-1})}}}
   {[m_{(1+B_{i-1})}]_q^{a_i}}\prod_{j=1+B_{i-1}}^{B_i}
   \frac{1}{[m_j]_q-\theta\,'
   q^{m_j}}\bigg\}\frac{q^{(a_s-1)m_{B_s}}}{[m_{B_s}]_q^{a_s-1}
   \big([m_{B_s}]_q-\theta\,'q^{m_{B_s}}\big)}\\
   &\quad-\theta\,'\sum_{m_1>\cdots>m_{B_s}>0}\bigg\{
   \prod_{i=1}^{s-1}\frac{q^{a_im_{(1+B_{i-1})}}}
   {[m_{(1+B_{i-1})}]_q^{a_i}}\prod_{j=1+B_{i-1}}^{B_i}
   \frac{1}{[m_j]_q-\theta\,'
   q^{m_j}}\bigg\}\frac{q^{a_sm_{B_s}}}{[m_{B_s}]_q^{a_s}
   \big([m_{B_s}]_q-\theta\,'q^{m_{B_s}}\big)}\\
   &\quad+\sum_{m_1>\cdots>m_{B_{s-1}}>0}\bigg\{
   \prod_{i=1}^{s-1}\frac{q^{a_im_{(1+B_{i-1})}}}
   {[m_{(1+B_{i-1})}]_q^{a_i}}\prod_{j=1+B_{i-1}}^{B_i}
   \frac{1}{[m_j]_q-\theta\,'
   q^{m_j}}\bigg\}\frac{q^{(a_s-1)m_{B_{s-1}}}}{[m_{B_{s-1}}]_q^{a_s}}\\
   &=f\bigg[\big(\prod_{i=1}^{s-1}x^{a_i}y^{b_i}\big)x^{a_s-1}y;\{0\}^s;
   \theta\,'\bigg]
   -\theta\,'f\bigg[\big(\prod_{i=1}^{s-1}x^{a_i}y^{b_i}\big)x^{a_s}y;
   \{0\}^s;\theta\,'\bigg]\\
   &\quad+f\bigg[\big(\prod_{i=1}^{s-2} x^{a_i}y^{b_i}\big)x^{a_{s-1}}
   y^{b_{s-1}-1}x^{a_s}y;\{0\}^s;\theta\,'\bigg]\\
   &\quad+(1-q)f\bigg[\big(\prod_{i=1}^{s-2} x^{a_i}y^{b_i}\big)x^{a_{s-1}}
   y^{b_{s-1}-1}x^{a_s-1}y;\{0\}^s;\theta\,'\bigg].
\end{align*}
\eop

\noindent{\bf Proof of Lemma~\ref{lem:14iiic}.} If $a>1$, then
\begin{align*}
    f\big[x^{a-1}y;\theta\big]-\theta f\big[x^ay;\theta\big]
   &= \sum_{m=1}^\infty
   \bigg(\frac{q^{(a-1)m}}{[m]_q^{a-1}}-\frac{\theta
   q^{am}}{[m]_q^a}\bigg)\frac{1}{[m]_q-\theta q^m}\\
   &=\sum_{m=1}^\infty \frac{q^{(a-1)m}}{[m]_q^a}\\
   &=\sum_{m=1}^\infty \frac{q^{(a-1)m}}{[m]_q^a}
   \cdot\frac{[m]_q-\theta\,'q^m}{[m]_q-\theta\,'q^m}\\
  &=\sum_{m=1}^\infty \frac{q^{(a-1)m}}{[m]_q^{a-1}\big([m]_q-
   \theta\,'q^m\big)}-\theta\,'\sum_{m=1}^\infty \frac{q^{am}}
   {[m]_q^a\big([m]_q-\theta\,'q^m\big)}\\
   &=f\big[x^{a-1}y;\theta\,'\big]-\theta\,'f\big[x^ay;\theta\,'\big].
\end{align*}
\eop

\noindent{\bf Proof of Lemma~\ref{lem:14iiid}.}  Let $n$ be a
positive integer.  Then
\begin{align*}
   &\theta\sum_{m=1}^n \frac{q^m}{[m]_q\big([m]_q-\theta q^m\big)}
   -\theta\,'\sum_{m=1}^n \frac{q^m}{[m]_q\big([m]_q-\theta
   q^m\big)}\\
   &= \sum_{m=1}^n
   \frac{[m]_q-\theta\,'q^m}{[m]_q\big([m]_q-\theta\,'q^m\big)}
   -\sum_{m=1}^n \frac{1}{[m]_q-\theta\,'q^m}-\sum_{m=1}^n
   \frac{[m]_q-\theta q^m}{[m]_q\big([m]_q-\theta q^m\big)}
   +\sum_{m=1}^n \frac{1}{[m]_q-\theta q^m}\\
   &=\sum_{m=1}^n \frac{1}{[m]_q-\theta q^m} -\sum_{m=1}^n
   \frac{1}{[m]_q-\theta\,'q^m}\\
   &=\sum_{m=0}^{n-1}\frac{1}{[m+1]_q-\theta q^{m+1}}-\sum_{m=1}^n
   \frac{1}{[m]_q-\theta\,'q^m}\\
   &= \sum_{m=0}^{n-1} \frac{1}{[m]_q-\theta\,'q^m}-\sum_{m=1}^n
   \frac{1}{[m]_q-\theta\,'q^m}\\
   &=\frac{1}{-\theta\,'} - \frac{1}{[n]_q-\theta\,'q^n}.
\end{align*}
The result now follows on letting $n\to\infty$. \eop

\section{Derivations}\label{sect:qDerivations}
We continue to employ the algebraic notation of the previous
section, and in addition, define the $\Q$-linear map
$\zeta^*=\lim_{q\to 1}\widehat{\zeta}$, so that
$\zeta^*(x^{s_1-1}y\cdots x^{s_m-1}y)=\zeta(s_1,\dots,s_m)$ gives
the ordinary multiple zeta value.  Note that $q$-duality
(Corollary~\ref{cor:qduality}) simply says that
$\widehat{\zeta}[\tau w]=\widehat{\zeta}[w]$ for all words
$w\in\mathfrak{h}^0$, while ordinary duality reduces to
$\zeta^*(\tau w)=\zeta^*(w)$. If $D$ is a derivation of
$\mathfrak{h}$, let $\overline{D}$ denote the conjugate derivation
$\tau D\tau$. As in~\cite{HoffOhno}, we refer to $D$ as symmetric
(resp.\ antisymmetric) if $\overline{D}=D$ $(\overline{D}=-D)$,
and note that any symmetric or antisymmetric derivation is
completely determined by where it sends $x$. Ihara and
Kaneko~\cite{Ihara} defined a family of antisymmetric derivations
$\partial_n$ for positive integers $n$ by declaring that
$\partial_n(x)=x(x+y)^{n-1}y$.  They conjectured---and
subsequently proved---that for all positive integers $n$ and words
$w\in\mathfrak{h}^0$, $\zeta^*(\partial_n(w))=0$.  Here, we shall
prove that this result extends to the multiple $q$-zeta function.

\begin{Thm}\label{thm:qDerivation} For all positive integers $n$ and words
$w\in\mathfrak{h}^0$, $\widehat{\zeta}[\partial_n(w)]=0$.
\end{Thm}

\noindent{\bf Proof.}  Again, for positive integer $n$ let $D_n$
be the derivation mapping $x\mapsto 0$ and $y\mapsto x^ny$.  Fix a
formal power series parameter $t$ and set
\[
   D:=\sum_{n=1}^\infty t^n\frac{D_n}{n},\qquad
   \sigma := \exp(D), \qquad
   \partial := \sum_{n=1}^\infty t^n\frac{\partial_n}{n}.
\]
The reformulated version of the generalized $q$-duality theorem
(Theorem~\ref{thm:f=g}) states that $\widehat{\zeta}[\sigma w] =
\widehat{\zeta}[\sigma\tau w]$ for all $w\in\mathfrak{h}^0$.  In
view of the special case, $q$-duality
(Corollary~\ref{cor:qduality}), this is equivalent to
$(\sigma-\overline{\sigma})w\in\mathrm{ker}\,\widehat{\zeta}$ for
all $w\in\mathfrak{h}^0$. We show that in fact,
$(\sigma-\overline{\sigma})\mathfrak{h}^0=\partial\mathfrak{h}^0$,
from which it follows that Theorem~\ref{thm:qDerivation} is
equivalent to generalized $q$-duality.    To prove the
equivalence, we require the following identity of Ihara and
Kaneko~\cite{Ihara}.

\begin{Prop}[Theorem 5.9 of~\cite{HoffOhno}]\label{prop:exp(partial)}
We have the following equality of $\mathfrak{h}[[t]]$
automorphisms: $\exp(\partial)=\overline{\sigma}\sigma^{-1}$.
\end{Prop}

To complete the proof of Theorem~\ref{thm:qDerivation}, observe as
in~\cite{Ihara,Zud} that since
\[
   \partial = \log\big(\overline{\sigma}\sigma^{-1}\big)
   =\log\big(1-(\sigma-\overline{\sigma})\sigma^{-1}\big)
   =-(\sigma-\overline{\sigma})\sum_{n=1}^\infty \frac1n
   \big((\sigma-\overline{\sigma})\sigma^{-1}\big)^{n-1}\sigma^{-1},
\]
and
\[
   \sigma-\overline{\sigma} =
   \big(1-\overline{\sigma}\sigma^{-1}\big)\sigma
   =\big(1-\exp(\partial)\big)\sigma
   =-\partial\sum_{n=1}^\infty \frac{\partial^{n-1}}{n!}\sigma,
\]
we see that $\partial\mathfrak{h}^0\subseteq
(\sigma-\overline{\sigma})\mathfrak{h}^0$ and
$(\sigma-\overline{\sigma})\mathfrak{h}^0\subseteq\partial\mathfrak{h}^0$.
Thus for the kernel of $\widehat{\zeta}$, we have the equivalences
\[
  (\sigma-\overline{\sigma})w\in \mathrm{ker}\,\widehat{\zeta}
  \quad
  \Longleftrightarrow
  \quad
  \partial w\in\mathrm{ker}\,\widehat{\zeta}
  \quad
  \Longleftrightarrow
  \quad  \forall n\in\Z^{+}, \;
  \widehat{\zeta}[\partial_nw]=0.
\]
\eop

\begin{Rem} The proof of
Proposition~\ref{prop:exp(partial)} given in~\cite{HoffOhno}
involves imposing a Hopf algebra structure on $\mathfrak{h}$ and
defining an action on it.  Zudilin~\cite[Lemma 7]{Zud} presents an
alternative proof in the case $t=1$ along the lines originally
indicated by Ihara and Kaneko~\cite{Ihara}.  It is possible to
extend Zudilin's presentation~\cite{Zud} to arbitrary $t$ by
defining a family $\{\varphi_s: s\in\R\}$ of automorphisms of
$\R\langle\langle x,y\rangle\rangle$ defined on the generators
$z=x+y$ and $y$ by
\[
   \varphi_s(z)=z,\qquad
   \varphi_s(y)=(1-tz)^sy\bigg(1-\frac{1-(1-tz)^s}{z}y\bigg)^{-1}.
\]
Routine calculations on the generators verify the equalities
\[
   \varphi_{s_1}\circ\varphi_{s_2}=\varphi_{s_1+s_2},
   \qquad \varphi_0=\mathrm{id},
   \qquad \frac{d}{ds}\varphi_s\bigg|_{s=0} = \partial,
   \qquad \varphi_1 = \overline{\sigma}\sigma^{-1}.
\]
The first three results imply that $\varphi_s = \exp(s\partial)$,
and the substitution $s=1$ gives
Proposition~\ref{prop:exp(partial)}.
\end{Rem}

\begin{Rem} In view of the identity $\partial_1=\overline{D}_1-D_1$,
the case $n=1$ of Theorem~\ref{thm:qDerivation} yields the
following $q$-analog of Hoffman's derivation theorem~\cite[Theorem
5.1]{Hoff92}~\cite[Theorem 2.1]{HoffOhno}:
\begin{Cor} For any word $w\in\mathfrak{h}^0$,
$\widehat{\zeta}[D_1w] = \widehat{\zeta}[\overline{D}_1 w]$.
Equivalently, if $s_1,\dots,s_m$ are positive integers with
$s_1>1$, then
\[
   \sum_{k=1}^m \zeta\big[\Cat_{j=1}^{k-1}s_j,1+s_k,\Cat_{j=k+1}^m
   s_j\big]
   = \sum_{k=1}^m \sum_{j=0}^{s_k-2}\zeta\big[\Cat_{i=1}^{k-1}s_i,
   s_k-j,j+1,\Cat_{i=k+1}^m s_i\big].
\]
By the usual convention on empty sums, the sum on the right is
zero if $s_k<2$.
\end{Cor}
\end{Rem}

\section{Cyclic Sums}\label{sect:cyclic}
In this section, we state and prove a $q$-analog of the cyclic sum
theorem~\cite{HoffOhno}, originally conjectured by Hoffman and
subsequently proved by Ohno using a partial fractions argument. As
a corollary, we give another proof of the $q$-sum formula
(Corollary~\ref{cor:qsum}).
\begin{Thm}[$q$-cyclic sum formula]\label{thm:cyclic} Let $n$ and $s_1,s_2,\dots,s_n$
be positive integers such that $s_j>1$ for some $j$.  Then
\[
   \sum_{j=1}^n \zeta\big[s_j+1,\Cat_{m=j+1}^n
   s_m,\Cat_{m=1}^{j-1}s_m\big]
   = \sum_{j=1}^n \sum_{k=0}^{s_j-2} \zeta\big[s_j-k,\Cat_{m=j+1}^n
   s_m,\Cat_{m=1}^{j-1} s_m,k+1\big].
\]
\end{Thm}
Note that the inner sum on the right vanishes if $s_j=1$.  We
refer to Theorem~\ref{thm:cyclic} as the $q$-cyclic sum formula
because, as with the limiting case in~\cite{HoffOhno}, it has an
elegant reformulation in terms of cyclic permutations of dual
argument lists.

\begin{Def} If $\vec s = (s_1,\dots,s_n)$ is a vector of $n$
positive integers, let
\[
   \mathscr{C}(\vec s) = \{(s_1,\dots,s_n),
   (s_2,\dots,s_n,s_1),\dots,(s_n,s_1,\dots,s_{n-1})\}
\]
denote the set of cyclic permutations of $\vec s$.  Also, for
notational convenience, define $\zeta^*[s_1,\dots,s_n]:=
\zeta[s_1+1,s_2,\dots,s_n]$.
\end{Def}
We can now restate Theorem~\ref{thm:cyclic} as follows.

\begin{Thm}[$q$-analog of~\cite{HoffOhno}, eq.~(2)]\label{thm:cyclicb}
Let $\mathbf s$ and $\mathbf s'$ be dual argument lists. Then
\[
   \sum_{\mathbf p\in\mathscr{C}(\mathbf s)}\zeta^*[\mathbf p]
   = \sum_{\mathbf p\in\mathscr{C}(\mathbf s')}\zeta^*[\mathbf p].
\]
\end{Thm}

To prove the implication Theorem~\ref{thm:cyclic}
$\Longrightarrow$ Theorem~\ref{thm:cyclicb}, we borrow an argument
of Ohno for the $q=1$ case.  Let
\[
   \mathbf s = \big(\Cat_{j=1}^m\{a_j+2,\{1\}^{b_j}\}\big) =
   (s_1,\dots,s_n),
\]
where $a_j$ and $b_j$ are non-negative integers for $1\le j\le m$
and $n=m+b_1+\dots+b_m$.  The right hand side of
Theorem~\ref{thm:cyclic} is
\begin{multline*}
   \frac{|\mathscr{C}(\mathbf s)|}{n}
   \sum_{\mathbf p\in\mathscr{C}(\mathbf s)}
   \sum_{k=0}^{p_1-2}\zeta[p_1-k,p_2,\dots,p_n,k+1]\\
   =\frac{|\mathscr{C}(\mathbf s)|}{n}
   \sum_{(\mathbf c,\mathbf d)}
   \sum_{k=0}^{c_1}\zeta[c_1+2-k,\{1\}^{d_1},\Cat_{j=2}^m
   \{c_j+2,\{1\}^{d_j}\},k+1],
\end{multline*}
where the outer sum on the right is over all cyclic permutations
\[
   (\mathbf c,\mathbf d)=((c_1,d_1),\dots,(c_m,d_m))
\]
of the ordered sequence of ordered pairs
$((a_1,b_1),\dots,(a_m,b_m))$.  Invoking $q$-duality
(Corollary~\ref{cor:qduality}), we find that the right hand side
of Theorem~\ref{thm:cyclic} can now be expressed as
\begin{align*}
   &\frac{|\mathscr{C}(\mathbf s)|}{n}
   \bigg\{\sum_{(\mathbf c,\mathbf d)}\zeta\big[\Cat_{j=1}^m
   \{c_j+2,\{1\}^{d_j}\},1\big]\\
   &\qquad +
   \sum_{(\mathbf c,\mathbf d)}\sum_{k=1}^{c_1}
   \zeta[c_1+2-k,\{1\}^{d_1},\Cat_{j=2}^m\{c_j+2,\{1\}^{d_j}\},k+1]\bigg\}\\
   &=\frac{|\mathscr{C}(\mathbf s)|}{n}
   \bigg\{\sum_{(\mathbf c,\mathbf d)}\zeta[d_m+3,\{1\}^{c_m},
   \Cat_{j=2}^m\{d_{m-j+1}+2,\{1\}^{c_{m-j+1}}\}]\\
   &\qquad+\sum_{(\mathbf c,\mathbf d)}\sum_{k=1}^{c_1}
   \zeta[2,\{1\}^{k-1},\Cat_{j=1}^{m-1}\{d_{m-j+1}+2,\{1\}^{c_{m-j+1}}\},
   d_1+2,\{1\}^{c_1-k}]\bigg\}\\
   &=\sum_{\mathbf p\in\mathscr{C}(\mathbf s')}\zeta^*[\mathbf p].
\end{align*}
But the left hand side of Theorem~\ref{thm:cyclic} is
\[
   \sum_{j=1}^n \zeta^*\big[\Cat_{m=j}^n
   s_m,\Cat_{m=1}^{j-1}s_m\big]
   = \sum_{\mathbf p\in\mathscr{C}(\mathbf s)} \zeta^*[\mathbf p].
\]

We now proceed with the proof of Theorem~\ref{thm:cyclic}.  As we
shall see, much of the proof of the limiting case
in~\cite{HoffOhno} can be adapted to the present situation with
only minor modifications.  To this end, we introduce two auxiliary
$q$-series.
\begin{Def}
For positive integers $s_1,s_2,\dots,s_n$ and non-negative integer
$s_{n+1}$, let
\begin{equation}
\begin{split}
   T[s_1,\dots,s_n] &:= \sum_{k_1>\cdots>k_{n+1}\ge 0}\;
   \frac{q^{k_1-k_{n+1}}}{[k_1-k_{n+1}]_q}\prod_{j=1}^n
   \frac{q^{(s_j-1)k_j}}{[k_j]_q^{s_j}},\\
   S[s_1,\dots,s_{n+1}] &:= \sum_{k_1>\cdots>k_{n+1}>0}\;
   \frac{q^{k_1}}{[k_1-k_{n+1}]_q}\prod_{j=1}^{n+1}
   \frac{q^{(s_j-1)k_j}}{[k_j]_q^{s_j}}.
\end{split}
\label{STdef}
\end{equation}
\end{Def}
For the convergence of the $q$-series~\eqref{STdef}, we have the
following generalization of~\cite[Theorem 3.1]{HoffOhno}.

\begin{Thm}\label{thm:TSconv} $T[s_1,\dots,s_n]$ is finite if
there is an index $j$ with $s_j>1$; $S[s_1,\dots,s_{n+1}]$ is
finite if one of $s_1,\dots,s_n$ exceeds 1 or if $s_{n+1}>0$.
\end{Thm}

We defer the proof of Theorem~\ref{thm:TSconv} to the end of the
section in order to proceed more directly with the proof of
Theorem~\ref{thm:cyclic}.  The key result we need is a direct
generalization of the corresponding result in~\cite{HoffOhno}:
\begin{Thm}[$q$-analog of~\cite{HoffOhno}, Theorem 3.2]\label{thm:Tdifference} If $s_1,\dots,s_n$ are positive
integers with $s_j>1$ for some $j$, then
\[
   T[s_1,\dots,s_n]-T[s_2,\dots,s_n,s_1]=
   \zeta[s_1+1,s_2,\dots,s_n] -
   \sum_{k=0}^{s_1-2}\zeta[s_1-k,s_2,\dots,s_n,k+1],
\]
where the sum on the right vanishes if $s_1=1$.
\end{Thm}
The proof of Theorem~\ref{thm:cyclic} now follows immediately on
summing Theorem~\ref{thm:Tdifference} over all cyclic permutations
of the argument sequence $s_1,\dots,s_n$.

\noindent{\bf Proof of Theorem~\ref{thm:Tdifference}.}  Although
we provide details, the argument is quite similar to the
corresponding argument in~\cite{HoffOhno}.  One minor difference
is that $\lim_{N\to \infty}1/[N]_q = 1-q\ne 0$ if $q\ne 1$, which
affects the computations used to arrive at~\eqref{S1=T} below.
First,
\begin{align}
   S[s_1,\dots,s_n,0] &= \sum_{k_1>\cdots>k_{n+1}>0}\;
   \frac{q^{k_1-k_{n+1}}}{[k_1-k_{n+1}]_q}\prod_{j=1}^n
   \frac{q^{(s_j-1)k_j}}{[k_j]_q^{s_j}}\nonumber\\
   &= \sum_{k_1>\cdots>k_{n+1}\ge 0}\;
   \frac{q^{k_1-k_{n+1}}}{[k_1-k_{n+1}]_q}\prod_{j=1}^n
   \frac{q^{(s_j-1)k_j}}{[k_j]_q^{s_j}}-
   \sum_{k_1>\cdots>k_n>0}\;\frac{q^{k_1}}{[k_1]_q}
   \prod_{j=1}^n \frac{q^{(s_j-1)k_j}}{[k_j]_q^{s_j}}\nonumber\\
   &= T[s_1,\dots,s_n] - \zeta[s_1+1,s_2,\dots,s_n].
\label{S0=T-zeta}
\end{align}
Next, we apply the identity
\begin{equation}
   \frac{q^{k_1-k_{n+1}}}{[k_1-k_{n+1}]_q[k_1]_q}
   = \frac{1}{[k_{n+1}]_q}\bigg(\frac{1}{[k_1-k_{n+1}]_q}-
   \frac{1}{[k_1]_q}\bigg)
\label{qparfrac}
\end{equation}
to $S[s_1,\dots,s_{n+1}]$. This gives
\begin{multline*}
  \sum_{k_1>\cdots>k_{n+1}>0}\;
   \frac{q^{k_1-k_{n+1}}}{[k_1-k_{n+1}]_q[k_1]_q}\cdot
   \frac{q^{(s_1-1)k_1}}{[k_1]_q^{s_1-1}}\cdot
   \frac{q^{s_{n+1}k_{n+1}}}{[k_{n+1}]_q^{s_{n+1}}}
   \prod_{j=2}^n \frac{q^{(s_j-1)k_j}}{[k_j]_q^{s_j}}\\
   =\sum_{k_1>\cdots>k_{n+1}>0}\;\bigg(\frac{1}{[k_1-k_{n+1}]_q}
   -\frac{1}{[k_1]_q}\bigg)\frac{q^{(s_1-1)k_1}}{[k_1]_q^{s_1-1}}
   \cdot\frac{q^{s_{n+1}k_{n+1}}}{[k_{n+1}]_q^{1+s_{n+1}}}
   \prod_{j=2}^n \frac{q^{(s_j-1)k_j}}{[k_j]_q^{s_j}},
\end{multline*}
from which it follows that
\begin{equation}
  S[s_1,\dots,s_{n+1}]
  =S[s_1-1,s_2,\dots,s_n,1+s_{n+1}]-\zeta[s_1,\dots,s_n,1+s_{n+1}].
\label{S=S-zeta}
\end{equation}
Finally, applying~\eqref{qparfrac} to
$S[1,s_2,\dots,s_n,s_{n+1}-1]$ gives
\begin{align*}
   &\sum_{k_1>\cdots>k_{n+1}>0}\;
   \frac{q^{k_1-k_{n+1}}}{[k_1-k_{n+1}]_q[k_1]_q}\cdot
   \frac{q^{(s_{n+1}-1)k_{n+1}}}{[k_{n+1}]_q^{s_{n+1}-1}}
   \prod_{j=2}^{n}\frac{q^{(s_j-1)k_j}}{[k_j]_q^{s_j}}\\
   &=\sum_{k_1>\cdots>k_{n+1}>0}\;\bigg(\frac{1}{[k_1-k_{n+1}]_q}-
   \frac{1}{[k_1]_q}\bigg)
   \prod_{j=2}^{n+1}\frac{q^{(s_j-1)k_j}}{[k_j]_q^{s_j}}\\
   &= \sum_{k_2>\cdots>k_{n+1}>0}\;
   \prod_{j=2}^{n+1}\frac{q^{(s_j-1)k_j}}{[k_j]_q^{s_j}}
   \lim_{N\to\infty}\sum_{k_1=k_2+1}^N\bigg(\frac{1}{[k_1-k_{n+1}]_q}-
   \frac{1}{[k_1]_q}\bigg)\\
   &=\sum_{k_2>\cdots>k_{n+1}>0}\;
   \prod_{j=2}^{n+1}\frac{q^{(s_j-1)k_j}}{[k_j]_q^{s_j}}
   \lim_{N\to\infty}\sum_{m=0}^{k_{n+1}-1}\bigg(\frac{1}{[k_2-m]_q}
   -\frac{1}{[N-m]_q}\bigg)\\
   &=\sum_{k_2>\cdots>k_{n+1}>0}\;
   \prod_{j=2}^{n+1}\frac{q^{(s_j-1)k_j}}{[k_j]_q^{s_j}}
   \sum_{m=0}^{k_{n+1}-1}\bigg(\frac{1}{[k_2-m]_q}+q-1\bigg)\\
   &=\sum_{k_2>\cdots>k_{n+1}>m\ge 0}\;
   \frac{q^{k_2-m}}{[k_2-m]_q}
   \prod_{j=2}^{n+1}\frac{q^{(s_j-1)k_j}}{[k_j]_q^{s_j}}.
\end{align*}
It follows that
\begin{equation}
   S[1,s_2,\dots,s_n,s_{n+1}-1] = T[s_2,\dots,s_n,s_{n+1}].
\label{S1=T}
\end{equation}
Now let $0\le j\le s_1-2$, apply~\eqref{S=S-zeta} and sum on $j$.
This yields
\[
   \sum_{j=0}^{s_1-2}S[s_1-j,s_2,\dots,s_n,j]
   = \sum_{j=0}^{s_1-2}\big(S[s_1-j-1,s_2,\dots,s_n,j+1]
   -\zeta[s_1-j,s_2,\dots,s_n,j+1]\big),
\]
which telescopes, leaving
\[
   S[s_1,s_2,\dots,s_n,0] = S[1,s_2,\dots,s_n,s_1-1]
   - \sum_{j=0}^{s_1-2}\zeta[s_1-j,s_2,\dots,s_n,j+1].
\]
Now apply~\eqref{S0=T-zeta} and~\eqref{S1=T} to obtain
\[
   T[s_1,\dots,s_n]-\zeta[s_1+1,s_2,\dots,s_n]
  =T[s_2,\dots,s_n,s_1]
   - \sum_{j=0}^{s_1-2}\zeta[s_1-j,s_2,\dots,s_n,j+1].
\]
\eop

As Ohno observed, the sum formula~\cite{Granville} is an easy
consequence of~\cite[Theorem 3.2]{HoffOhno}.  Correspondingly, we
can give another proof of Corollary~\ref{cor:qsum}, our $q$-analog
of the sum formula.

\noindent{\bf Alternative Proof of Corollary~\ref{cor:qsum}.} Sum
Theorem~\ref{thm:Tdifference} over all $s_1,\dots,s_n$ with
$s_1+\cdots+s_n=k$.  Since the resulting sum of $T$-functions
vanishes, we get
\begin{align*}
   \sum_{s_1+\cdots+s_n=k}\zeta[s_1+1,s_2,\dots,s_n]
   &=
   \sum_{s_1+\cdots+s_n=k}\sum_{j=0}^{s_1-2}\zeta[s_1-j,s_2,\dots,s_n,
   j+1]\\
   &= \sum_{s_1+\cdots+s_{n+1}=k}\zeta[s_1+1,s_2,\dots,s_{n+1}].
\end{align*}
It follows that the sums are independent of $n$; whence each is
equal to
\[
   \sum_{s_1=k}\zeta[s_1+1] = \zeta[k+1],
\]
as required. \eop

We conclude the section with a proof of Theorem~\ref{thm:TSconv}.
Again, the argument closely follows Ohno's proof of the limiting
case in~\cite{HoffOhno}.

\noindent{\bf Proof of Theorem~\ref{thm:TSconv}.}
By~\eqref{S0=T-zeta},
\[
   S[s_1,\dots,s_n,s_{n+1}]\le S[s_1,\dots,s_n,0]\le
   T[s_1,\dots,s_n],
\]
so $S[s_1,\dots,s_{n+1}]$ is finite if $T[s_1,\dots,s_n]$ is.
By~\eqref{S1=T},
\[
   S[1,s_2,\dots,s_n,s_{n+1}]=T[s_2,\dots,s_n,s_{n+1}+1],
\]
so the statement about finiteness of $S$ follows from the
corresponding statement about $T$.  To prove finiteness of
$T[s_1,\dots,s_n]$ with $s_1+\cdots+s_n>n$, it suffices to
consider the case $s_1+\cdots+s_n=n+1$, for if $s_k>1$, then
$T[s_1,\dots,s_n]\le T[\{1\}^{k-1},2,\{1\}^{n-k}]$.  Thus, we need
only prove that $T[\{1\}^{k-1},2,\{1\}^{n-k}]<\infty$ for $1\le
k\le n$.  When $k=1$, we have
\begin{align*}
  T[2,\{1\}^{n-1}] &= \sum_{k_1>\cdots>k_{n+1}\ge 0}\;
  \frac{q^{k_1-k_{n+1}+k_1}}{[k_1-k_{n+1}]_q[k_1]_q^2}
  \prod_{j=2}^n \frac{1}{[k_j]_q}\\
  &\le \sum_{\substack{k_1>\cdots>k_n>0\\k_1\ge m>0}}
  \frac{q^{m+k_1}}{[m]_q[k_1]_q^2}\prod_{j=2}^n
  \frac{1}{[k_j]_q}\\
  &= \zeta[3,\{1\}^{n-1}]+n\zeta[2,\{1\}^n] +
  \sum_{k=1}^{n-1}\zeta[2,\{1\}^{k-1},2,\{1\}^{n-k-1}]\\
  &<\infty.
\end{align*}
Arguing inductively, we now suppose that
$T[\{1\}^{k-1},2,\{1\}^{n-k}]<\infty$ for some $k\ge 1$.
By~\eqref{S0=T-zeta},~\eqref{S1=T} and the inductive hypothesis,
\begin{align*}
   T[\{1\}^k,2,\{1\}^{n-k-1}] &= S[\{1\}^k,2,\{1\}^{n-k-1},0] +
   \zeta[2,\{1\}^{k-1},2,\{1\}^{n-k-1}]\\
   &=
   T[\{1\}^{k-1},2,\{1\}^{n-k}]+\zeta[2,\{1\}^{k-1},2,\{1\}^{n-k-1}]\\
   &<\infty,
\end{align*}
as required. \eop

\section{Multiple $q$-Polylogarithms}\label{sect:polylogs}
In analogy with~\cite[eq.~(1.1)]{BBBLa}, define
\begin{equation}
\label{qpolylog}
   \lambda_q\bigg[\begin{array}{cc}s_1,\dots,s_m \\ b_1,\dots,b_m \end{array}
   \bigg] :=
   \sum_{\nu_1,\dots,\nu_m>0}\;\prod_{k=1}^m
   b_k^{-\nu_k}\bigg[\sum_{j=k}^m \nu_j\bigg]_q^{-s_k},
\end{equation}
and set
\begin{equation}
\label{qLi}
   \mathrm{Li}_{s_1,\dots,s_m}[x_1,\dots,x_m] :=
   \sum_{n_1>\cdots>n_m>0}\;\prod_{k=1}^m
   \frac{x_k^{n_k}}{[n_k]_q^{s_k}}.
\end{equation}
The substitution $n_k=\sum_{j=k}^m \nu_j$ shows
that~\eqref{qpolylog} and~\eqref{qLi} are related by
\[
   \mathrm{Li}_{s_1,\dots,s_m}[x_1,\dots,x_m] =
   \lambda_q\bigg[\begin{array}{cc}s_1,\dots,s_m\\
   y_1,\dots,y_m\end{array}\bigg],
   \qquad y_k=\prod_{j=1}^k x_j^{-1}.
\]
\begin{Thm}[$q$-analog of Theorem 9.1 of~\cite{BBBLa}]\label{thm:multisect}
Let $b_1,\dots,b_m\in\C$, $s_1,\dots,s_m>0$ and let $n$ be a
positive integer. Then
\[
   n^m\lambda_{q^n}\bigg[\begin{array}{cc}s_1,\dots,s_m\\
   b_1^n,\dots,b_m^n\end{array}\bigg]
   = [n]_q^s\sum_{\varepsilon_1^n=\cdots=\varepsilon_m^n=1}
   \lambda_q\bigg[\begin{array}{cc}s_1,\dots,s_m\\ \varepsilon_1b_1,
   \dots,\varepsilon_mb_m\end{array}\bigg],
\]
where the sum is over all $n^m$ sequences
$(\varepsilon_1,\dots,\varepsilon_m)$ of complex $n^{\mathrm{th}}$
roots of unity, and $s=\sum_{k=1}^m s_k$.
\end{Thm}
\noindent{\bf Proof.}  In light of the identity
\[
   \frac{1}{[\nu]_{q^n}^s} =
   \bigg(\frac{1-q^n}{1-q^{n\nu}}\bigg)^s
   =
   \bigg(\frac{1-q^n}{1-q}\bigg)^s\bigg(\frac{1-q}{1-q^{n\nu}}\bigg)^s
   =\frac{[n]_q^s}{[n\nu]_q^s},
\]
we have
\begin{align*}
   n^m\lambda_{q^n}\bigg[\begin{array}{cc}s_1,\dots,s_m\\
   b_1^n,\dots,b_m^n\end{array}\bigg]
   &= n^m\sum_{\nu_1,\dots,\nu_m>0}\;\prod_{k=1}^m b_k^{-n\nu_k}
   \bigg[\sum_{j=k}^m \nu_j\bigg]_{q^n}^{-s_j}\\
   &= n^m\sum_{\nu_1,\dots,\nu_m>0}\;\prod_{k=1}^m b_k^{-n\nu_k}
   [n]_q^{s_j}\bigg[n\sum_{j=k}^m \nu_j\bigg]_q^{-s_j}\\
   &= [n]_q^{s}\sum_{\nu_1,\dots,\nu_m>0}\;\prod_{k=1}^m nb_k^{-n\nu_k}
   \bigg[\sum_{j=k}^m n\nu_j\bigg]_q^{-s_j}\\
   &= [n]_q^{s}\sum_{\nu_1,\dots,\nu_m>0}\;\prod_{k=1}^m
   b_k^{-\nu_k} \bigg[\sum_{j=k}^m \nu_j\bigg]_q^{-s_j}
   \sum_{\mu_k=0}^{n-1} e^{-2\pi i \mu_k\nu_k/n}\\
   &=[n]_q^{s}\sum_{\mu_1=0}^{n-1}\cdots\sum_{\mu_m=0}^{n-1}\;
   \sum_{\nu_1,\dots,\nu_m>0}\;\prod_{k=1}^m
   b_k^{-\nu_k} e^{-2\pi i \mu_k \nu_k/n}
   \bigg[\sum_{j=k}^m \nu_j\bigg]_q^{-s_j}.
\end{align*}
Letting $\varepsilon_k = e^{2\pi i \mu_k/n}$ completes the
proof.\eop

In contrast with our proof of Theorem~\ref{thm:multisect}, the
proof of the limiting case in~\cite{BBBLa} made use of the
Drinfel'd simplex integral representation for multiple
polylogarithms.  As integral representations for multiple
polylogarithms have proved eminently useful in establishing many
of their properties, we derive here a $q$-analog of the Drinfel'd
simplex integral for the multiple
$q$-polylogarithm~\eqref{qpolylog}.

\begin{Thm}\label{thm:qintegral}
Let $s_1,\dots,s_m$ be positive integers.  For the multiple
$q$-polylogarithm, we have the multiple Jackson $q$-integral
representation
\begin{equation}
\label{Jack}
   \lambda_q\bigg[\begin{array}{cc}s_1,\dots,s_m\\y_1,\dots,y_m\end{array}
   \bigg]
   = \int \prod_{k=1}^m \bigg(\prod_{r=1}^{s_k-1}\frac{d_q
   t_r^{(k)}}{t_r^{(k)}}\bigg)\frac{d_q
   t_{s_k}^{(k)}}{y_k-t_{s_k}},
\end{equation}
where the multiple Jackson $q$-integral~\eqref{Jack} is over the
simplex
\[
   1>t_1^{(1)}>\cdots>t_{s_1}^{(1)}>\cdots>t_1^{(m)}>\cdots>t_{s_m}^{(m)}>0.
\]
\end{Thm}
\begin{Rem} As in~\cite{BBBLa}, we may abbreviate~\eqref{Jack} by
\[
    \lambda_q\bigg[\begin{array}{cc}s_1,\dots,s_m\\y_1,\dots,y_m\end{array}
   \bigg]
   =(-1)^m \int_0^1 \prod_{k=1}^m (\omega[0])^{s_k-1}\omega[y_k],
   \qquad \omega[b]:=\frac{d_q t}{t-b}.
\]
\end{Rem}
\begin{Cor} For multiple $q$-zeta values, we have the multiple
Jackson $q$-integral representation
\[
   \zeta[s_1,\dots,s_m] = (-1)^m \int_0^1 \prod_{k=1}^m
   (\omega[0])^{s_k-1}\omega\bigg[\prod_{j=1}^k q^{1-s_j}\bigg].
\]
\end{Cor}

\noindent{\bf Proof of Theorem~\ref{thm:qintegral}.} We first
establish the following
\begin{Lem}  Let $s$ be a positive integer, $0<t_0<1$ and $m>0$.
Then
\[
   \int\limits_{t_0>t_1>\cdots>t_s>0}
   \bigg(\prod_{r=1}^{s-1}\frac{d_q t_r}{t_r}\bigg)t_s^{m-1}d_q
   t_s = \frac{t_0^m}{[m]_q^s}.
\]
\end{Lem}
\noindent{\bf Proof.}  When $s=1$, the integral reduces to the
geometric series
\[
   \int\limits_{t_0>t_1>0} t_1^{m-1}d_q t_1 =
   (1-q)t_0\sum_{j=0}^\infty q^j (q^j t_0)^{m-1} =
   \bigg(\frac{1-q}{1-q^m}\bigg)t_0^m.
\]
Suppose the lemma holds for $s-1$.  By the inductive hypothesis,
\[
   \int\limits_{t_0>t_1>\cdots>t_s>0}
   \bigg(\prod_{r=1}^{s-1}\frac{d_q t_r}{t_r}\bigg)t_s^{m-1}d_q
   t_s = \int\limits_{t_0>t_1>0}\frac{t_1^m}{[m]_q^{s-1}}
   \frac{d_q t_1}{t_1}
   = \frac{1}{[m]_q^{s-1}}\int\limits_{t_0>t_1>0} t_1^{m-1}d_q t_1
   = \frac{t_0^m}{[m]_q^s},
\]
as required. \eop

To prove~\eqref{Jack}, it will suffice to
establish the identity
\begin{equation}
\label{doit}
   \int \prod_{k=1}^m \bigg(\prod_{r=1}^{s_k-1}\frac{d_q
   t_r^{(k)}}{t_r^{(k)}}\bigg)\frac{d_q
   t_{s_k}^{(k)}}{y_k-t_{s_k}^{(k)}}
   = \lambda_q\bigg[\begin{array}{cc}s_1,\dots,s_m\\
   y_1/t_0,\dots,y_m/t_0\end{array}\bigg],
\end{equation}
where the integral~\eqref{doit} is over the simplex
\[
t_0>t_1^{(1)}>\cdots>t_{s_1}^{(1)}>\cdots>t_1^{(m)}>\cdots>t_{s_m}^{(m)}>0.
\]
When $m=1$,~\eqref{doit} reduces to
\begin{align*}
   \int\limits_{t_0>t_1>\cdots>t_s>0}\bigg(\prod_{r=1}^{s-1}
   \frac{d_q t_r}{t_r}\bigg)\frac{y^{-1}d_q t_s}{1-y^{-1}t_s}
   &=\int\limits_{t_0>t_1>\cdots>t_s>0}\bigg(\prod_{r=1}^{s-1}
   \frac{d_q t_r}{t_r}\bigg)\sum_{\nu=1}^\infty
   y^{-\nu}t_s^{\nu-1}d_q t_s\\
   &= \sum_{\nu=1}^\infty y^{-\nu}\int\limits_{t_0>\cdots>t_s>0}
   \bigg(\prod_{r=1}^{s-1}\frac{d_q t_r}{t_r}\bigg)t_s^{\nu-1}d_q
   t_s\\
   &= \sum_{\nu=1}^\infty \frac{y^{-\nu}t_0^{\nu}}{[\nu]_q^s}\\
   &= \lambda\bigg[\begin{array}{c}s\\y/t_0\end{array}\bigg].
\end{align*}
Suppose~\eqref{doit} holds for $m-1$.  Then the inductive
hypothesis implies that the integral~\eqref{doit} is equal to
\begin{align*}
   &\int\limits_{t_0>t_1>\cdots>
   t_{s_1}>0}\bigg(\prod_{r=1}^{s_1-1} \frac{d_q t_r}{t_r}\bigg)
   \frac{y_1^{-1}d_q t_{s_1}}{1-y_1^{-1}t_{s_1}}
   \sum_{\nu_2,\dots,\nu_m>0}\;
   \prod_{k=2}^m t_{s_1}^{\nu_k} y_k^{-\nu_k}
   \bigg[\sum_{j=k}^m \nu_j\bigg]_q^{-s_j}\\
   &= \sum_{\nu_1,\dots,\nu_m>0}y_1^{-\nu_1}\prod_{k=2}^m
   y_k^{-\nu_k}\bigg[\sum_{j=k}^m \nu_j\bigg]_q^{-s_j}
   \int\limits_{t_0>t_1>\cdots>t_{s_1}>0}
   \bigg(\prod_{r=1}^{s_1-1}\frac{d_q t_r}{t_r}\bigg)
   t_{s_1}^{\nu_1+\nu_2+\cdots+\nu_m-1}d_q t_{s_1}\\
   &= \sum_{\nu_1,\dots,\nu_m>0}\prod_{k=1}^m y_k^{-\nu_k}
   t_0^{\nu_k}\bigg[\sum_{j=k}^m \nu_j\bigg]_q^{-s_k}\\
   &= \lambda_q\bigg[\begin{array}{cc}s_1,\dots,s_m\\
   y_1/t_0,\dots,y_m/t_0\end{array}\bigg].
\end{align*}
\eop

\begin{Rem} Zhao~\cite{Zhao} has outlined an alternative approach
to deriving the multiple Jackson $q$-integral representation of
the multiple $q$-polylogarithm.  In addition, he studies the
$q$-shuffles, first explicated in~\cite[Section 7]{BowBradSurvey},
that arise when multiplying two such integrals.
\end{Rem}

\section{A Double Generating Function for $\zeta[m+2,\{1\}^n]$}
\label{sect:qEuler}
In this section, we derive the
following $q$-analog of~\cite[eq.~(10)]{BBB} and a few of its
implications.

\begin{Thm}\label{thm:Drin} The double generating function
identity
\begin{multline}\label{Drin}
   \sum_{m=0}^\infty\sum_{n=0}^\infty u^{m+1}v^{n+1}
   \zeta[m+2,\{1\}^n]\\ = 1-\exp\bigg\{\sum_{k=2}^\infty
   \big\{u^k+v^k-\big(u+v+(1-q)uv\big)^k\big\}\frac1k\sum_{j=2}^k
   (q-1)^{k-j}\zeta[j]\bigg\}
\end{multline}
holds.
\end{Thm}
Noting that the generating function~\eqref{Drin} is symmetric in
$u$ and $v$, we immediately derive the following special case of
$q$-duality.

\begin{Cor} For all non-negative integers $m$ and $n$,
$\zeta[m+2,\{1\}^n] = \zeta[n+2,\{1\}^m]$.
\end{Cor}

Of course, we have already proved $q$-duality at full strength
(Corollary~\ref{cor:qduality}) as a consequence of generalized
$q$-duality (Theorem~\ref{thm:genqduality}).  The main interest
for Theorem~\ref{thm:Drin} may be that it shows that
$\zeta[m+2,\{1\}^n]$ can be expressed in terms of sums of products
of depth-1 $q$-zeta values. When $n=1$, this reduces to the
following convolution identity, which provides a $q$-analog of
Euler's evaluation~\cite[eq.~(31)]{BBB}~\cite{LE,Niels} of
$\zeta(m+2,1)$.

\begin{Cor}\label{cor:qEuler} Let $m$ be a non-negative integer.  Then
\[
   2\zeta[m+2,1] = (m+2)\zeta[m+2] +(1-q)\,m\,\zeta[m+2] -
   \sum_{k=2}^{m+1}\zeta[m+3-k]\,\zeta[k].
\]
\end{Cor}

In particular, when $m=0$ we get $\zeta[2,1]=\zeta[3]$, which
corrects an error in~\cite[Theorem 15]{Zud}.

\noindent{\bf Proof of Corollary~\ref{cor:qEuler}.}  Compare
coefficients of $u^{m+1}v^2$ on each side of the double generating
function identity~\eqref{Drin}.  Letting
\[
   c_k := \begin{cases}\displaystyle
   \frac1k\sum_{j=2}^k(q-1)^{k-j}\zeta[j], &\text{if $k\ge 2$}\\
   0, &\text{if $k<2$,}
   \end{cases}
\]
we find that
\begin{multline}\label{extract}
   2\zeta[m+2,1] =
   (m+2)c_{m+3}+2(1-q)(m+1)c_{m+2}+(1-q)^2\,m\,c_{m+1}\\
   -\sum_{k+l=m+3}c_kc_l+2(q-1)\sum_{k+l=m+2}c_kc_l
   -(q-1)^2\sum_{k+l=m+1}c_kc_l,
\end{multline}
where convolution sums in~\eqref{extract} range over all integers
$k$ and $l$ satisfying the indicated relations.  Now
\begin{align*}
   &(m+2)c_{m+3}+2(1-q)(m+1)c_{m+2}+(1-q)^2\,m\,c_{m+1}\\
   &= (m+2)\sum_{j=2}^{m+3}(q-1)^{m+3-j}\zeta[j]-2(m+1)
   \sum_{j=2}^{m+2}(q-1)^{m+3-j}\zeta[j]
   +m\sum_{j=2}^{m+1}(q-1)^{m+3-j}\zeta[j]\\
   &=
   \big\{(m+2)-2(m+1)+m\big\}\sum_{j=2}^{m+1}(q-1)^{m+3-j}\zeta[j]\\
   &\qquad\qquad\qquad+(m+2)\sum_{j=m+2}^{m+3}(q-1)^{m+3-j}\zeta[j]
   -2(m+1)(q-1)\zeta[m+2]\\
   &=(m+2)\zeta[m+3]+(1-q)m\zeta[m+2].
\end{align*}
In light of~\eqref{extract}, it now follows that
\begin{multline*}
   2\zeta[m+2,1] - (m+2)\zeta[m+3] - (1-q)m\zeta[m+2]\\
   = -\sum_{k+l=m+3}c_kc_l+2(q-1)\sum_{k+l=m+2}c_kc_l
   -(q-1)^2\sum_{k+l=m+1}c_kc_l.
\end{multline*}
To avoid having to deal directly with boundary cases, we set
$\zeta_{+}[n] := \zeta[n](n\ge 2)$ and $(q-1)_{+}^n = (q-1)^n(n\ge
0)$.   Then
\begin{align*}
   &2\zeta[m+2,1] - (m+2)\zeta[m+3] - (1-q)m\zeta[m+2]\\
  &=-\sum_{k\in\Z} c_{m+3-k}
     \bigg\{c_k-2(q-1)c_{k-1}+(q-1)^2c_{k-2}\bigg\}\\
  &=-\sum_{k\in\Z}c_{m+3-k}\sum_{j\in\Z}\bigg\{(q-1)_{+}^{k-j}
  -2(q-1)_{+}^{k-1-j}+(q-1)^2(q-1)_{+}^{k-2-j}\bigg\}\zeta[j]\\
  &=-\sum_{k\in\Z}c_{m+3-k}\bigg\{\zeta_{+}[k]+
  \big\{(q-1)-2(q-1)\big\}\zeta_{+}[k-1]\\
     &\qquad\qquad\qquad+\sum_{j\le k-2}
     \big\{(q-1)^{k-j}-2(q-1)^{k-j}+(q-1)^{k-j}\big\}\zeta_{+}[j]
     \bigg\}\\
  &=\sum_{k\in\Z}c_{m+3-k}(q-1)\zeta_{+}[k-1]
  -\sum_{k\in\Z}c_{m+3-k}\,\zeta_{+}[k].
\end{align*}
We now re-index the latter two sums, replacing $k$ by $m+4-n$ in
the first, and $k$ by $m+3-n$ in the second.  Thus,
\begin{align*}
   & 2\zeta[m+2,1] - (m+2)\zeta[m+3] - (1-q)m\zeta[m+2]\\
   &=\sum_{n\in\Z}\zeta_{+}[m+3-n](q-1)c_{n-1}-\sum_{n\in\Z}
  \zeta_{+}[m+3-n]c_n\\
  &=\sum_{n\in\Z}\zeta_{+}[m+3-n]\sum_{j\in\Z}\big\{(q-1)(q-1)_{+}^{n-1-j}
  -(q-1)_{+}^{n-j}\big\}\zeta_{+}[j]\\
  &=\sum_{n\in\Z}\zeta_{+}[m+3-n]\bigg\{\sum_{j\le
  n-1}\big\{(q-1)^{n-j}-(q-1)^{n-j}\big\}\zeta_{+}[j]-(q-1)_{+}^0
  \zeta_{+}[n]\bigg\}\\
  &= -\sum_{n\in\Z}\zeta_{+}[m+3-n]\zeta_{+}[n]\\
  &=-\sum_{n=2}^{m+1}\zeta[m+3-n]\zeta[n],
\end{align*}
as claimed. \eop

\begin{Rem} Similarly, one could derive an explicit identity for
$\zeta[m+2,1,1]$ in terms of depth-1 $q$-zeta values by comparing
coefficients of $u^{m+1}v^3$ in Theorem~\ref{thm:Drin}.  The
resulting identity would be a $q$-analog of Markett's double
convolution identity~\cite{Markett} for $\zeta(m+2,1,1)$.
\end{Rem}

Our proof of Theorem~\ref{thm:Drin} employs techniques from the
theory of basic hypergeometric series.  For real $x$ and $y$ and
non-negative integer $n$, the asymmetric $q$-power~\cite{QCalc} is
given by
\[
   (x+y)_q^n := \prod_{k=0}^{n-1} (x+yq^k),
   \qquad (x+y)_q^\infty := \lim_{n\to\infty}(x+y)_q^n.
\]
The $q$-gamma function~\cite[p.~493]{AndAskRoy}~\cite[p.~16]{Gasp}
is defined by
\[
   \Gamma_q(x)=\frac{(1-q)_q^\infty
   (1-q)^{1-x}}{(1-q^x)_q^\infty},
\]
and the basic hypergeometric
function~\cite[p.~520]{AndAskRoy}~\cite[p.\ xv, eq.\ (22)]{Gasp}
is
\[
   {}_2\phi_1\bigg[\begin{array}{cc}q^a,
   q^b\\q^c\end{array}\bigg|x\bigg]
   = \sum_{n=0}^\infty \frac{(1-q^a)_q^n
   (1-q^b)_q^n}{(1-q^c)_q^n (1-q)_q^n}\,x^n,
   \qquad |x|<1.
\]
Heine's $q$-analog of Gauss's summation formula for the ordinary
hypergeometric function~\cite[p.~522]{AndAskRoy}~\cite[p.\ xv,
eq.\ (23)]{Gasp} may be stated in the form
\begin{equation}\label{Heine}
   {}_2\phi_1\bigg[\begin{array}{cc}q^a,
   q^b\\q^c\end{array}\bigg|q^{c-a-b}\bigg]
   =
   \frac{\Gamma_q(c)\Gamma_q(c-a-b)}{\Gamma_q(c-a)\Gamma_q(c-b)},
   \qquad \big|q^{c-a-b}\big|<1.
\end{equation}

Our first step towards proving Theorem~\ref{thm:Drin} is to
establish the following result.

\begin{Thm}[$q$-analog of eq.~(6.5) of~\cite{BBBLa}]\label{thm:2phi1}
Let $x$ and $y$ be real numbers satisfying $|x|<1$ and $|y|<1$.
Then
\begin{equation}
\label{2phi1}
   \sum_{m=0}^\infty \sum_{n=0}^\infty (-1)^{m+n} [x]_q^{m+1}
   [y]_q^{n+1} \zeta[m+2,\{1\}^n] =
   1-\frac{\Gamma_q(1+x)\Gamma_q(1+y)}{\Gamma_q(1+x+y)}.
\end{equation}
\end{Thm}

\noindent{\bf Proof of Theorem~\ref{thm:2phi1}.}  Let $L$ denote
the bivariate double generating function on the left hand side
of~\eqref{2phi1}.  Then
\begin{align*}
    L &=-[y]_q\sum_{m=0}^\infty
   (-1)^{m+1}[x]_q^{m+1}\sum_{k=1}^\infty
   \frac{q^{(m+1)k}}{[k]_q^{m+2}}\prod_{j=1}^{k-1}\bigg(1-\frac{[y]_q}{[j]_q}
   \bigg)\\
   &=-[y]_q\sum_{m=0}^\infty
  (-1)^{m+1}[x]_q^{m+1}\sum_{k=1}^\infty
   \frac{q^{(m+1)k}}{[k]_q^{m+2}}\prod_{j=1}^{k-1}\frac{[j]_q-[y]_q}{[j]_q}\\
   &=-[y]_q\sum_{m=0}^\infty
   (-1)^{m+1}[x]_q^{m+1}\sum_{k=1}^\infty
   \frac{q^{(m+1)k}}{[k]_q^{m+2}}\prod_{j=1}^{k-1}\frac{q^y-q^j}{1-q^j}\\
   &= q^y\bigg(\frac{1-q^{-y}}{1-q}\bigg)\sum_{m=0}^\infty
   (-1)^{m+1}[x]_q^{m+1}\sum_{k=1}^\infty
   \frac{q^{(m+1)k}}{[k]_q^{m+2}}\cdot\frac{q^{(k-1)y}}{(1-q)_q^{k-1}}
   \prod_{j=1}^{k-1}\big(1-q^{j-y}\big)\\
      &=\sum_{m=0}^\infty
   (-1)^{m+1}[x]_q^{m+1}\sum_{k=1}^\infty
  \frac{q^{(m+1)k}q^{ky}}{[k]_q^{m+1}}\cdot\frac{(1-q^{-y})_q^k}{(1-q)_q^k}.
\end{align*}
Now interchange order of summation, noting that the sum on $m$ is
a geometric series.  Thus, we find that
\begin{align*}
  L &= \sum_{k=1}^\infty \frac{q^{ky}(1-q^{-y})_q^k}{(1-q)_q^k}
  \sum_{m=0}^\infty
  \frac{(-1)^{m+1}q^{(m+1)k}[x]_q^{m+1}}{[k]_q^{m+1}}\\
  &=-\sum_{k=1}^\infty \frac{q^{ky}(1-q^{-y})_q^k}{(1-q)_q^k}\cdot
  \frac{q^k[x]_q/[k]_q}{1+q^k[x]_q/[k]_q}\\
 &= -\sum_{k=1}^\infty \frac{q^{(y+1)k}(1-q^{-y})_q^k}{(1-q)_q^k}
  \cdot\frac{[x]_q}{[k]_q+q^k[x]_q}\\
    &= -\sum_{k=1}^\infty
  \frac{q^{(y+1)k}(1-q^{-y})_q^k}{(1-q)_q^k}\cdot
  \frac{1-q^x}{1-q^{k+x}}\\
  &= -\sum_{k=1}^\infty\frac{q^{(y+1)k}(1-q^{-y})_q^k}{(1-q)_q^k}\cdot
  \frac{(1-q^x)_q^k}{(1-q^{1+x})_q^k}\\
  &=
  1-{}_2\phi_1\bigg[\begin{array}{cc}q^{-y},q^x\\q^{1+x}
  \end{array}\bigg|q^{1+y}\bigg].\\
\end{align*}
Invoking Heine's formula~\eqref{Heine} completes the proof. \eop

To express the right hand side of~\eqref{2phi1} in the form of an
exponentiated power series, we require the following series
expansion of the logarithm of the $q$-gamma function.

\begin{Lem}\label{lem:logqGamma} For real $x$ such that $-1<x<1$,
we have
\[
   \log \Gamma_q(1+x) = -\gamma_q\, x + \sum_{k=2}^\infty
    \frac{[x]_q^k}{k}\sum_{j=2}^k (q-1)^{k-j}\zeta[j],
\]
where
\[
   \gamma_q := \log(1-q)-\frac{\log q}{1-q}\sum_{n=1}^\infty
   \frac{q^n}{[n]_q}
\]
is a $q$-analog of Euler's constant, $\gamma$.
\end{Lem}

\noindent{\bf Proof of Lemma~\ref{lem:logqGamma}.}  By definition,
\[
   \Gamma_q(1+x)=(1-q)^{-x}\prod_{n=1}^\infty
   \frac{1-q^n}{1-q^{n+x}}.
\]
Therefore,
\begin{align}\label{qGammatilde}
   \log\Gamma_q(1+x)+x\log(1-q)
   &= -\sum_{n=1}^\infty \log\bigg(\frac{1-q^{n+x}}{1-q^n}\bigg)
   = -\sum_{n=1}^\infty
   \log\bigg(1+\bigg(\frac{1-q^x}{1-q^n}\bigg)q^n\bigg)\nonumber\\
   &=\sum_{n=1}^\infty \sum_{k=1}^\infty
   \frac{(-1)^k[x]_q^k\,q^{kn}}{k\, [n]_q^k} \nonumber\\
   &=\sum_{k=1}^\infty \frac{(-1)^k [x]_q^k}{k}\,\widetilde{\zeta}[k],
\end{align}
where
\[
   \widetilde{\zeta}[k]:= \sum_{n=1}^\infty \frac{q^{kn}}{[n]_q^k},
   \qquad k>0.
\]
If we now multiply the identity
\[
    \frac{q^{jn}}{[n]_q^j} = (q-1)\frac{q^{(j-1)n}}{[n]_q^{j-1}} +
    \frac{q^{(j-1)n}}{[n]_q^j}, \qquad  (n,j\in\Z^{+})
\]
by $(q-1)^{k-j}$ and sum on $n$ and $j$, we find that
\[
   \sum_{j=2}^k (q-1)^{k-j}\,\widetilde{\zeta}[j] =
   \sum_{j=2}^k (q-1)^{k-j+1}\,\widetilde{\zeta}[j-1]
   +\sum_{j=2}^k (q-1)^{k-j}\zeta[j],
\]
which telescopes, leaving us with
\begin{equation}\label{tilde}
   \widetilde{\zeta}[k] = (q-1)^{k-1}\,\widetilde{\zeta}[1]
   +\sum_{j=2}^k (q-1)^{k-j}\zeta[j],
   \qquad k\ge 1.
\end{equation}
If we now substitute~\eqref{tilde} into~\eqref{qGammatilde}, there
comes
\begin{align*}
   &\log\Gamma_q(1+x)+x\log(1-q)\\
   &= \sum_{k=1}^\infty \frac{(-1)^k[x]_q^k}{k}\bigg\{
   (q-1)^{k-1}\,\widetilde{\zeta}[1]+\sum_{j=2}^k
   (q-1)^{k-j}\zeta[j]\bigg\}\\
   &= -(q-1)^{-1}\,\widetilde{\zeta}[1]\log(1+(q-1)[x]_q)
   +\sum_{k=1}^\infty
   \frac{(-1)^k[x]_q^k}{k}\sum_{j=2}^k(q-1)^{k-j}\zeta[j]\\
   &= \frac{x\widetilde{\zeta}[1]\log q}{1-q}+\sum_{k=2}^\infty
   \frac{(-1)^k[x]_q^k}{k}\sum_{j=2}^k(q-1)^{k-j}\zeta[j].
\end{align*}
In light of the fact that
\[
   \log\Gamma_q(1+x)=-\gamma x+\sum_{k=2}^\infty \frac{(-1)^k
   x^k}{k}\zeta(k)
\]
and $\lim_{q\to 1-}\Gamma_q(1+x)=\Gamma(1+x)$, it follows that
$\lim_{q\to 1-} \gamma_q= \gamma$.  Thus, the proof of
Lemma~\ref{lem:logqGamma} is complete. \eop

\noindent{\bf Proof of Theorem~\ref{thm:Drin}.}  By
Theorem~\ref{thm:2phi1} and Lemma~\ref{lem:logqGamma}, we have
\begin{multline*}
   \sum_{m=0}^\infty\sum_{n=0}^\infty (-1)^{m+n}[x]_q^{m+1}
   [y]_q^{n+1}\zeta[m+2,\{1\}^n]\\
   = 1-\exp\bigg\{\sum_{k=2}^\infty\frac{(-1)^k}{k}
   \big([x]_q^k+[y]_q^k-[x+y]_q^k\big)\sum_{j=2}^k
   (q-1)^{k-j}\zeta[j]\bigg\}.
\end{multline*}
Noting that $[x+y]_q = [x]_q+[y]_q+(q-1)[x]_q[y]_q$, the result
now follows on replacing $[x]_q$ by $-u$ and $[y]_q$ by $-v$. \eop

\section{Acknowledgment} The author is grateful to Masanobu
Kaneko, Yasuo Ohno, and Jun-ichi Okuda for making preprints of
their work available to him, and for their kindness and
hospitality.  Thanks are due also to Mike Hoffman and David
Broadhurst for several engaging discussions.

\end{document}